\documentclass[11pt,a4paper]{article}

\usepackage{amsfonts,amssymb,amsmath}
\usepackage{graphicx}
\usepackage{epstopdf}
\usepackage{booktabs}
\usepackage{makecell}
\usepackage{braket}
\usepackage{bm}
\usepackage{algorithmic}
\usepackage{xcolor,colortbl}
\usepackage[hypertexnames=false,colorlinks=true,linkcolor=blue,citecolor=blue,
            bookmarksopen=false,bookmarks=false,
            pdfstartview=XYZ,pdffitwindow=true,pdfcenterwindow=true]{hyperref}
\ifpdf
  \DeclareGraphicsExtensions{.eps,.pdf,.png,.jpg}
\else
  \DeclareGraphicsExtensions{.eps}
\fi
\usepackage{pifont}
\newcommand{\cmark}{\ding{51}}

\newcommand{\eps}{\varepsilon}
\newcommand{\mdsr}[1]{\textcolor{black}{#1}}

\title{Towards a new classification of bursting patterns:\\ review \& extensions}

 \author{
   Mathieu Desroches\thanks{MathNeuro Team, Inria Sophia Antipolis M{\'e}diterran{\'e}e Research Centre, 2004 route des lucioles BP93, 06902 Sophia Antipolis cedex, France and MCEN Team, Basque Centre for Applied Mathematics (BCAM), Mazarredo Kalea, 21, Bilbao 48009, Bizkaia, Spain, \href{mathieu.desroches@inria.fr}{mathieu.desroches@inria.fr}}
   \and
   John Rinzel\thanks{Center for Neural Science, New York University, 4 Washington Pl, New York, NY 10003, USA and Courant Institute for Mathematical Sciences, New York University, 251 Mercer St, New York 10012, USA}
   \and
   Serafim Rodrigues\thanks{Ikerbasque, MCEN Team, Basque Centre for Applied Mathematics (BCAM), Alameda de Mazarredo, 14, Bilbao 48009, Bizkaia, Spain, \href{srodrigues@bcamath.org}{srodrigues@bcamath.org}}
}

\usepackage{amsopn}
\usepackage{mathtools}

\definecolor{Darkgray}{rgb}{0.4,0.4,0.4}

\newcommand{\blue}[1]{\textcolor{blue}{#1}}
\newcommand{\red}[1]{\textcolor{black}{#1}}
\newcommand{\bluebis}[1]{\textcolor{black}{#1}}

\begin{document}
\maketitle

\begin{abstract}
The mathematical classification of complex bursting oscillations in multiscale excitable systems, seen for example in physics and neuroscience, has been the subject of active enquiry since the early 1980s and is still ongoing. This classification problem is fundamental as it also establishes analytical and numerical foundations for studying complex temporal behaviours in multiple timescale models. This manuscript begins by reviewing the seminal works of Rinzel~\cite{rinzel86} and Izhikevich~\cite{izhikevich00} in classifying bursting patterns of excitable cell models. Moreover, we recall an alternative, yet complementary, mathematical classification approach by Golubitsky, which together with the Rinzel-Izhikevich proposals provide the state-of-the-art foundations to the classification problem. Unexpectedly, while keeping within the Rinzel-Izhikevich framework, we find novel cases of bursting mechanisms, which were not considered before. Moving beyond the state-of-the-art, we identify novel bursting mechanisms that fall outside the Rinzel, Izhikevich and Golubitsky classification system. This leads us towards a new classification, which requires the analysis of both the fast and the slow subsystems of an underlying slow-fast model. This new classification allows the dynamical dissection of a larger class of bursters. To substantiate this, we add a new class of bursters with at least two slow variables, which we denote \textit{folded-node bursters}, to convey the idea that the bursts are initiated or annihilated via a \textit{folded-node} singularity. In fact, there are two main families of folded-node bursters, depending upon the phase of the bursting cycle during which folded-node dynamics occurs. If it occurs during the silent phase, we obtain the \textit{classical} folded-node bursting (or simply folded-node bursting). If it occurs during the active phase, we have \textit{cyclic} folded-node bursting. We classify both families and give examples of minimal systems displaying these novel types of bursting behaviour. 
\end{abstract}

\section{Introduction}
\label{sec:intro}
The fascination of experimentalists, physicists and mathematicians towards spontaneous and complex oscillations dates back to the early nineteenth century, particularly through observations of electro-chemical systems~\cite{HedgesMyers1926}. Indeed, how can seemingly ``inert sub-components'' assemble into ``life'', in what is currently understood (in biophysics) as open multi-scale (far from equilibrium) systems with dissipative structures. Van der pol was among the first scientists to exhibit equations with multiple timescales and a dissipative structure, which display oscillations akin to those observed in electro-chemical systems and that indeed could not be explained by previous mathematical theories~\cite{appleton22,vanderpol26,vanderpol28}. Despite remarkable advances, it is only relatively recently (since the 1980s) that a deeper understanding of nonlinear multi-scale complex oscillations was made possible due to the development of a coherent mathematical theories and classification system. The present manuscript follows up on our recent work on multi-scale systems~\cite{desroches12,desroches16b}, further builds upon several seminal works and finally proposes an extended classification framework, which we envisage will guide future developments of analytical, numerical and modelling work on multi-scale systems.\newline

To contain the complexity of mathematical characterisation of multi-scale systems, the present manuscript will solely focus on the dynamics emerging from the interplay between the timescales of a system modelled by a slow-fast ODE (Ordinary Differential Equations). While this approach circumvents the wider unresolved mathematical barriers in generally describing multi-scale systems across spatial-temporal scales (e.g. via partial differential equations), it will significantly enable us to obtain a deeper insight on emergent timescale-induced dynamics. This will later inform these other multi-scale approaches. Moreover, despite the relative simplicity of slow-fast ODEs, the associated theory is still in development. The present manuscript aims to contribute to it, even in a small way. More importantly, slow-fast ODEs have already enabled remarkable predictions of complex oscillations, hence their relevance. Specifically, these slow-fast ODEs are written in a singularly perturbed ODE form as follows:
\begin{equation}\label{eq:canburstsl}
\begin{aligned}
\eps\dot{x} &= f(x,y)\\
~\dot{y} &= g(x,y),
\end{aligned}
\end{equation}
where $(x, y) \in \mathbb{R}^n \times \mathbb{R}^m$, $f$ and $g$ are smooth functions, the dot denotes differentiation with respect to the slow time $t$ and $0\!<\!\eps\!\ll\!1$ is a small time constant that explicitly distinguishes the different timescales between the variables $x$ and $y$. Here $x$ are the fast processes and $y$ the slow ones. Indeed, this modelling approach captures many complex phenomena (e.g. biological, physical, chemical) where the timing and ordering of its underlying processes ensures correct functions. For example, in enzyme kinetics the binding of the substrate to the active site on the enzyme occurs within a few nano-seconds, while enzymatic reactions last up to a few seconds. This mix of various timescales appears to give rise to sudden, often surprising, jumps in the state of biological systems and, in general, to complex oscillations. An alternative way of writing equation~\eqref{eq:canburstsl} is to note that, as long as $\eps\neq 0$, one can obtain an equivalent system by rescaling time to the fast time $\tau=t/\eps$, which gives the form
\begin{equation}\label{eq:canburstfa}
\begin{aligned}
x' &= f(x,y)\\
y' &= \eps g(x,y),
\end{aligned}
\end{equation}
where the prime denotes differentiation with respect to $\tau$. A priori, these equivalences do not reveal any new information since we observe the same system but with a different ``magnifying glass'' (i.e. a different scaling). However, although the systems~\eqref{eq:canburstsl} and~\eqref{eq:canburstfa}, remain equivalent for every $\eps\!>\!0$, they have different set of mathematical equations in the (singular) limit $\eps=0$ (detailed below). This highlights the fact that these systems are \textit{singularly} and not \textit{regularly} perturbed. Resolving these singular limits and reconciling these dual mathematical equations has been the basis of several mathematical and numerical developments within, e.g., \textit{Catastrophe theory}, \textit{Singularity theory}, \textit{Bifurcation theory}, \textit{Geometrical Singular perturbation theory (Fenichel theory)}, \textit{geometric desingularization or blow-up method}. In short, these two different limits give rise to the so-called \textit{slow subsystem} and \textit{fast subsystem}, respectively. The slow subsystem is a differential-algebraic problem, where the slow variables remain explicitly dynamic through the differential equation $\dot{y}=g(x,y)$, while the fast variables are enslaved to the slow ones through the algebraic equation $f(x,y)=0$. On the other hand, the fast subsystem is a family of dynamical systems on the fast variables $x$, where the slow variables $y$ have lost their dynamics and have become parameters. The set $\{f=0\}$ is referred to as the \textit{critical manifold} of the system, and it plays a central role in both subsystems. Specifically, it is the phase space of the slow subsystem and it is the set of equilibria of the fast subsystem. Hence the locus of the limiting slow motion corresponds to a subset of the bifurcation diagram of the fast subsystem obtained when varying one or more slow variables. Moreover, varying other system parameters can induce the limiting slow motion to have non trivial trajectories, such as \textit{canard solutions} that emerge via \textit{dynamic bifurcations} and that visit both the stable and unstable regions of the critical manifold (see e.g.~\cite{desroches12} for details). 
\begin{figure}[!t]
\centering
\includegraphics{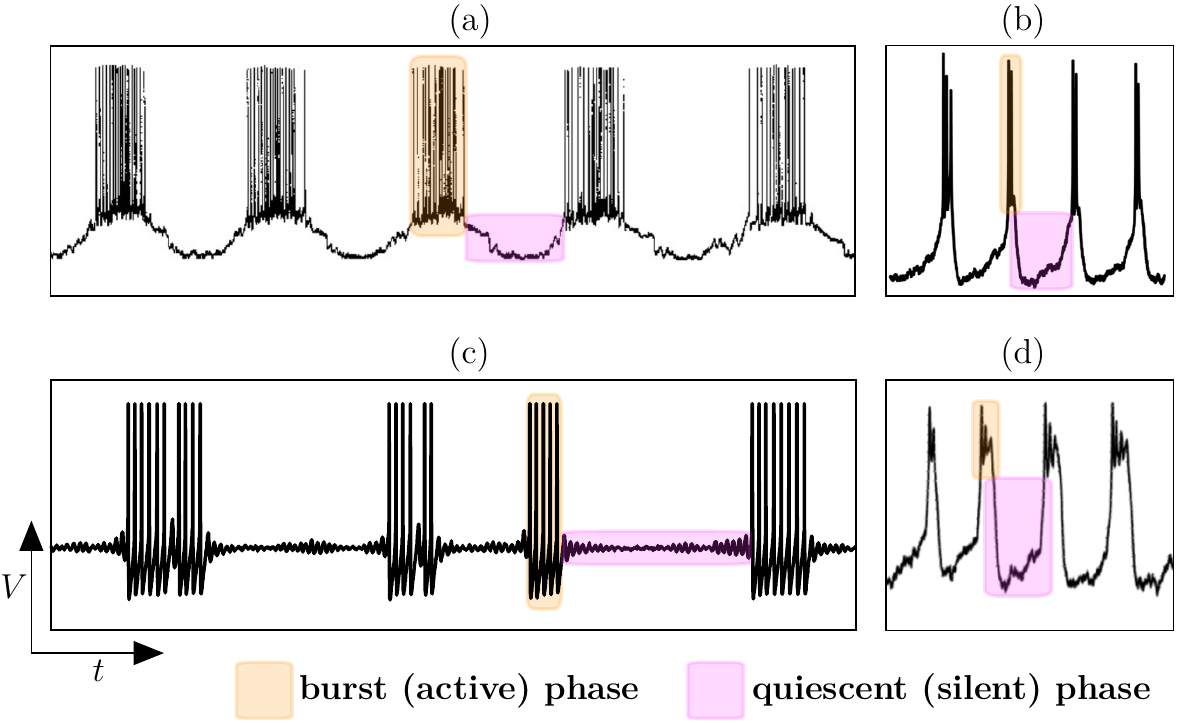}
\caption{Example of electrophysiological recordings of bursting oscillations in four types of neurons: (a) parabolic-type bursting from the CeN neuron from the melibe (a sea slug)~cite{newcomb08}; (b) square-wave-type bursting from a human $\beta$-cell~\cite{riz14}; (c) elliptic-type bursting from a dorsal-root-ganglia (DRG) neuron of a rat~\cite{jian04}; (d) Pseudo-plateau-type bursting from a pituitary cell of a rat~\cite{tabak11}.}
\label{fig.bursting_data}
\end{figure}
The knowledge of this dissection between the slow and fast subsystems, together with the knowledge of how the trajectory of the full system (i.e. for small $\eps\!>\!0$) evolves along these subsystems is at the basis of the various classification systems for complex slow-fast oscillations that we will subsequently review. Specifically, we will focus on the state-of-the-art mathematical classification systems for \textit{bursting dynamics}, their limitations and finally propose a novel classification framework. \red{Our novel classification fundamentally relies upon a certain type of canard configuration, referred to as \textit{folded node}. Noteworthy, canards are indirectly included in the previous classification frameworks as boundaries between the spiking and the bursting regimes via so-called \textit{spike-adding transitions} or \textit{torus canards}, which we will review, thus highlighting the importance of canards in the classification process. We will define the notion of bursting oscillations in the context of neuronal systems} (see examples in Figure~\ref{fig.bursting_data}) since for historical reasons the notion of bursting emerged within the neuroscientific literature. In particular, \mdsr{bursting models} appeared in the context of classical single neuron electrophysiological measurements, where \mdsr{the neuron's voltage time-series \red{displays} a bursting oscillation either in response to a brief input stimulus or, in absence of any stimuli, in an endogenous manner}. These oscillations are defined as having a periodic succession (sometimes irregular) of two distinct epochs of activity. One epoch features slow and low-amplitude activity, and it is typically referred to as the \textit{quiescent} (or \textit{silent}) phase. The other epoch features fast and high-amplitude activations (i.e. several action potentials or spikes) is classically denoted \textit{active} or \textit{burst} phase as shown by several examples in Figure~\ref{fig.bursting_data}. Although a great deal of our discussions will be in the context of neuronal dynamics, the mathematical framework intends to capture complex slow-fast oscillations beyond the scope of neuroscientific applications (e.g. in chemical reactions, genetic switches, material transitions, etc.). Hence some of the mathematical model constructs that we will present here display bursting oscillations not necessarily observed in neural data. Indeed, our \mdsr{idealized} models will not have direct biophysical interpretation as we aim to be as general as possible in describing the fundamental mathematical mechanisms, which can then be applied to explain complex bursting oscillations in multiple contexts. Moreover, we will focus on the minimal deterministic mathematical setting for bursting oscillations. Specifically, we will cover the case of two-timescale systems with explicit timescale separation (dictated by a single small parameter $0\!<\!\eps\!\ll\!1$), with two fast variables that will enable the description of the active phase of a bursting oscillation (i.e. $x\in\mathbb{R}^2$), and one or two slow variables describing the quiescent phase of the bursting dynamics (i.e. $y\in\mathbb{R}$ or $y\in\mathbb{R}^2$). This minimal setting will inform more complex scenarios involving multi-dimensional systems with multiple timescales.

\mdsr{This paper is organized as follows. In Section~\ref{sec:review}, we will review existing classification frameworks for bursting oscillations. \red{Subsequently}, in Section~\ref{sec:beyond}, we first introduce the \red{key} idea of our novel bursting classification based upon the concept of folded-node bursting dynamics. \red{This is followed by showcasing} several \red{new} examples of folded-node burster idealized models, first in the case of classical folded node and then in the case of cyclic folded node. Finally, in the conclusion section, we review our findings and propose a number of perspectives and future directions to explore. We \red{complete the manuscript} by proposing, in Appendix~\ref{sec:newfast}, few \red{novel} additional bursting scenarios within the Rinzel-Izhikevich's classification. \red{These include cases} with transcritical and pitchfork bifurcations of limit cycles, \textit{isola bursting} and a two-slow-variable bursting scenario with a family of transcritical bifurcations of equilibria.}

\begin{figure}[!t]
\centering
\includegraphics{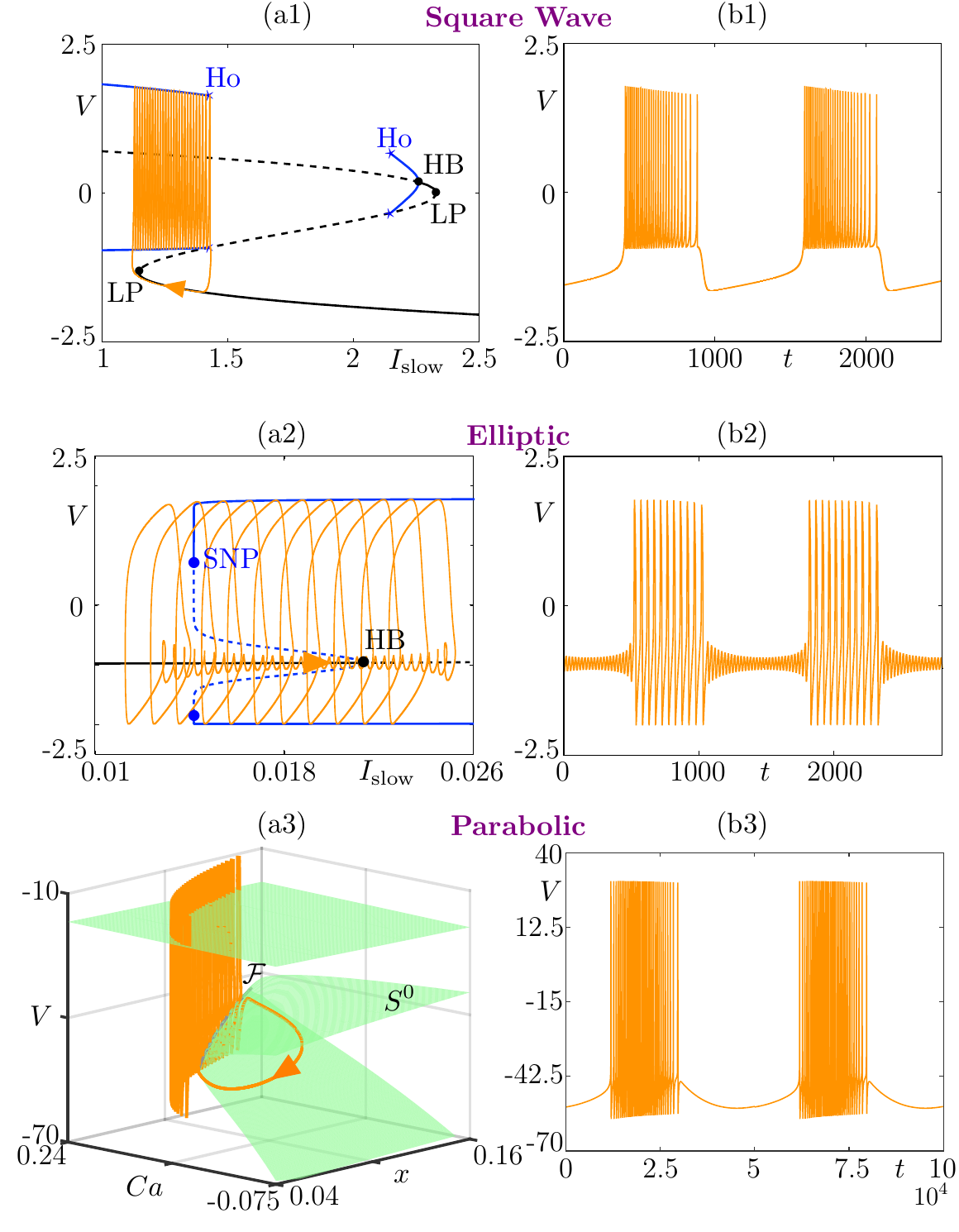}
\caption{Rinzel classification of bursting patterns: square-wave bursting, here in the Hindmarsh-Rose model~\cite{hindmarsh84} (panels (a1)-(b1)); elliptic bursting, here in the FitzHugh-Rinzel model~\cite{rinzel86,rinzel87} (panels (a2)-(b2)); parabolic bursting, here in Plant's model~\cite{plant81} (panels (a3)-(b3)).}
\label{fig:rinzel}
\end{figure}

\section{Review of the state-of-the-art classification of bursting patterns}
\label{sec:review}
%
\subsection{Rinzel's classification (mid 1980s)} Historically, John Rinzel opened the door towards mathematically understanding bursting oscillations. His seminal work on a mathematical analysis and classification of bursting oscillatory patterns, were first published within two companion manuscripts~\cite{rinzel86,rinzel87}. The fundamental insight behind Rinzel's classification is based on slow-fast dissection and in particular describing the bifurcation structure of the fast subsystem where the slow variables are frozen. Subsequently, the time trajectory of the full system (i.e. for small $\eps\!>\!0$) is superimposed on top of the bifurcation structure of the fast subsystem. This reveals that the quiescent phase of the bursting cycle correspond to trajectory segments where the solution slowly tracks families of stable equilibria, or low-amplitude (subthreshold) limit cycles, of the fast subsystem. Conversely, the burst phase of the full system's cycle correspond to trajectory segments where the solution slowly tracks families of limit cycles of the fast subsystem. The transitions between these two main phases of bursting cycles occur near bifurcation points of the fast subsystem. With this approach, Rinzel proposed three classes of bursting dynamics based on both the bifurcation structure of the fast subsystem and the salient features of the main fast variable's time profile (in the neuronal context this is typically the neuronal membrane potential). These features include spike frequency during the burst, dynamics during the silent phase (oscillatory or not), shape of the burst (on a plateau compared to the silent phase or on the contrary with undershoots). These three features led Rinzel to name three classes as \textit{square-wave}, \textit{elliptic} and \textit{parabolic} bursting. We show an example of each class in Figure~\ref{fig:rinzel}.

\subsection{Izhikevich's classification (ca. 2000)} 
\begin{figure}[!h]
\centering
\includegraphics{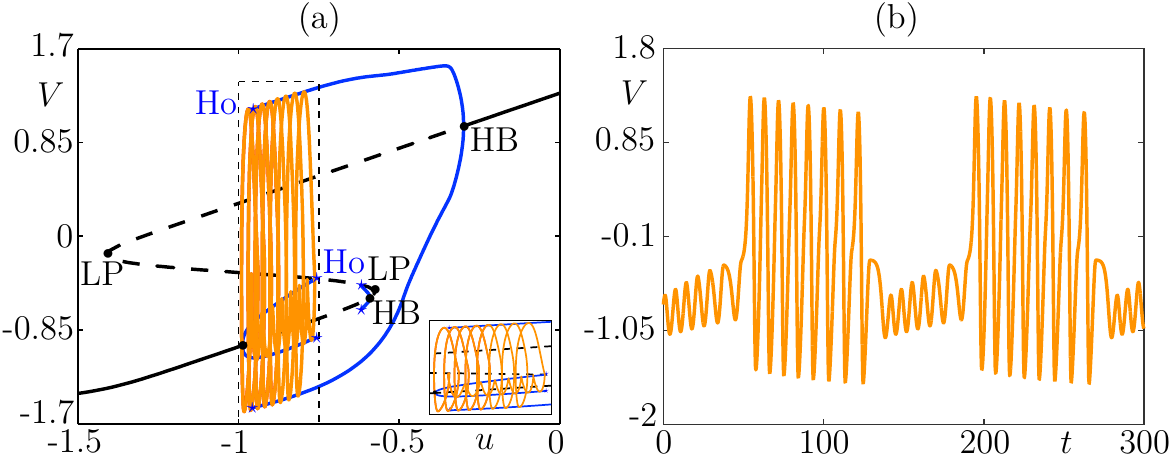}
\caption{Small homoclinic / big homoclinic bursting, corresponding to Fig.88 of~\cite{izhikevich00}. Panel (a) shows the slow-fast dissection in the $(u,V)$ phase plane and panel (b) shows the $V$-time series of this bursting solution.}
\label{fig:homhom}
\end{figure}

Eugene Izhikevich generalised Rinzel's approach by considering that a bursting pattern is entirely characterised by a pair of bifurcations ($\mathbf{b_1}$, $\mathbf{b_2}$) of the fast subsystem. One bifurcation, say $\mathbf{b_1}$, explains the transition from quiescence to burst, and the other, $\mathbf{b_2}$, marks the inverse transition, from burst to quiescence. Due to the well established bifurcation theory and indeed knowledge of classes of bifurcation, this led to a systematic identification of at least 120 bursting patterns~\cite{izhikevich00}. An example of a bursting model that is not within the Rinzel classification is depicted in Fig~\ref{fig:homhom}. In this example the bursting pattern has a transition from quiescence to burst via a homoclinic bifurcation (involving a small homoclinic connection) and equally, the transition from burst to quiescence is via homoclinic bifurcation (involving a large homoclinic connection). In many ways, Izhikevich's work serves as a key source of reference for classification of complex slow-fast oscillations. This is particularly the case in Neuroscience since some of the assembled examples were motivated by existing conductance-based neuronal models and demonstrated how complex neuronal oscillations could be achieved by adding one slow equation to a spiking system. Indeed, a dedicated book towards Neuroscience was later published, where the derived models where also put into context with neurophysiological processes~\cite{izhikevich07}. The result of this deeply insightful work is a quasi-complete classification of bursting patterns in terms of pairs of codimension-one bifurcations of the fast subsystem.
%
\subsection{Bertram et al.'s / Golubitsky et al.'s classification (mid 1990)} An alternative approach to classification was proposed by Bertram and colleagues in 1995~\cite{bertram95} and extended mathematically by Golubitsky and collaborators in 2001~\cite{golubitsky01} using a singularity theory viewpoint. The fundamental idea consists in identifying a codimension-$k$ bifurcation point ($k \geq 2$) in the fast subsystem and subsequently consider the slow variables of the bursting system as the unfolding parameters of this high codimension bifurcation point. The bursting is then obtained via a slow path made by the slow variables in the unfolding of that point (i.e. within a multidimensional parameter space). The minimum codimension, whose unfolding allows to create a given bursting pattern, defines the class of the associated bursting patterns provided a notion of path equivalence is properly defined. Specifically, two paths are equivalent if one can pass from one to the other via diffeomorphism and a re-parameterization. Recently, a review and a show-case demonstrating the construction of bursting oscillations via this approach, including cases for higher codimension bifurcation points was published in~\cite{saggio17}. It is worth noting that the Rinzel-Izhikevich approach and the Bertram-Golubitsky approach both focus on the fast subsystem only. Moreover, a way to see a link between the two approaches is that the two bifurcation points ($\mathbf{b_1}$, $\mathbf{b_2}$) of the fast subsystem (as characterised by Izhikevich's approach) belong to bifurcation curves in a two-parameter plane, which coalesce at a codimension-two bifurcation point that characterises this particular bursting pattern from the singularity theory viewpoint. This implies that in principle the Rinzel-Izhikevich and the Bertram-Golubitsky approaches both lead to the same number of bursting oscillation cases.

\mdsr{The bursting patterns covered by these three existing classification schemes have not been exhausted yet, even though a large number (way above one hundred) have already been reported and analysed in previous works. However, we propose a few more cases which we believe have not been considered before and which will be presented with associated idealized models in Appendix~\ref{sec:newfast}. In particular, we will show bursting scenarios where the burst phase ends due to a transcritical or a pitchfork bifurcation of limit cycles of the fast subsystem. We also propose the concept of \textit{isola bursting}, where the burst starts and ends through saddle-node bifurcations of limit cycles of the fast subsystem which happen to lie on an isola of limit cycles. Finally, we propose one example (amongst many) of bursting pattern with two slow variables where the burst initiates through a family of transcritical bifurcation of equilibria.}

\section{\mdsr{Towards a new classification of bursting patterns}}
\label{sec:beyond}

\subsection{\mdsr{Main idea to go beyond the state-of-the-art}}
\begin{figure}[!t]
\centering
\includegraphics{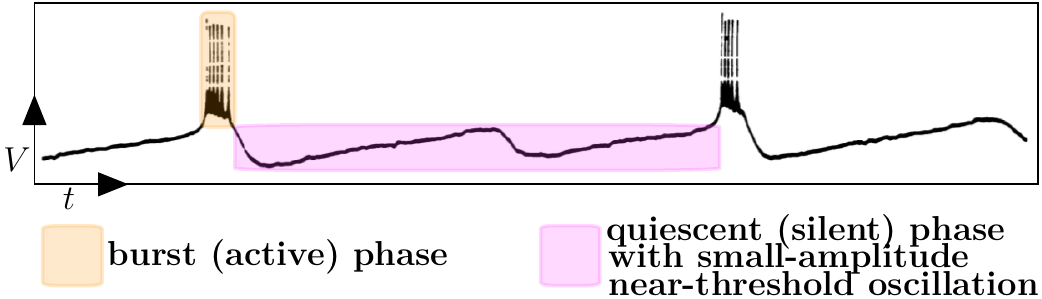}
\caption{Electrophysiological recordings of the lateral thalamic nuclei neuron in cat from~\cite{roy84} show complex bursting oscillations.}
\label{fig:MMBO_exp}
\end{figure}

\begin{figure}[!t]
\centering
\includegraphics[width=14cm]{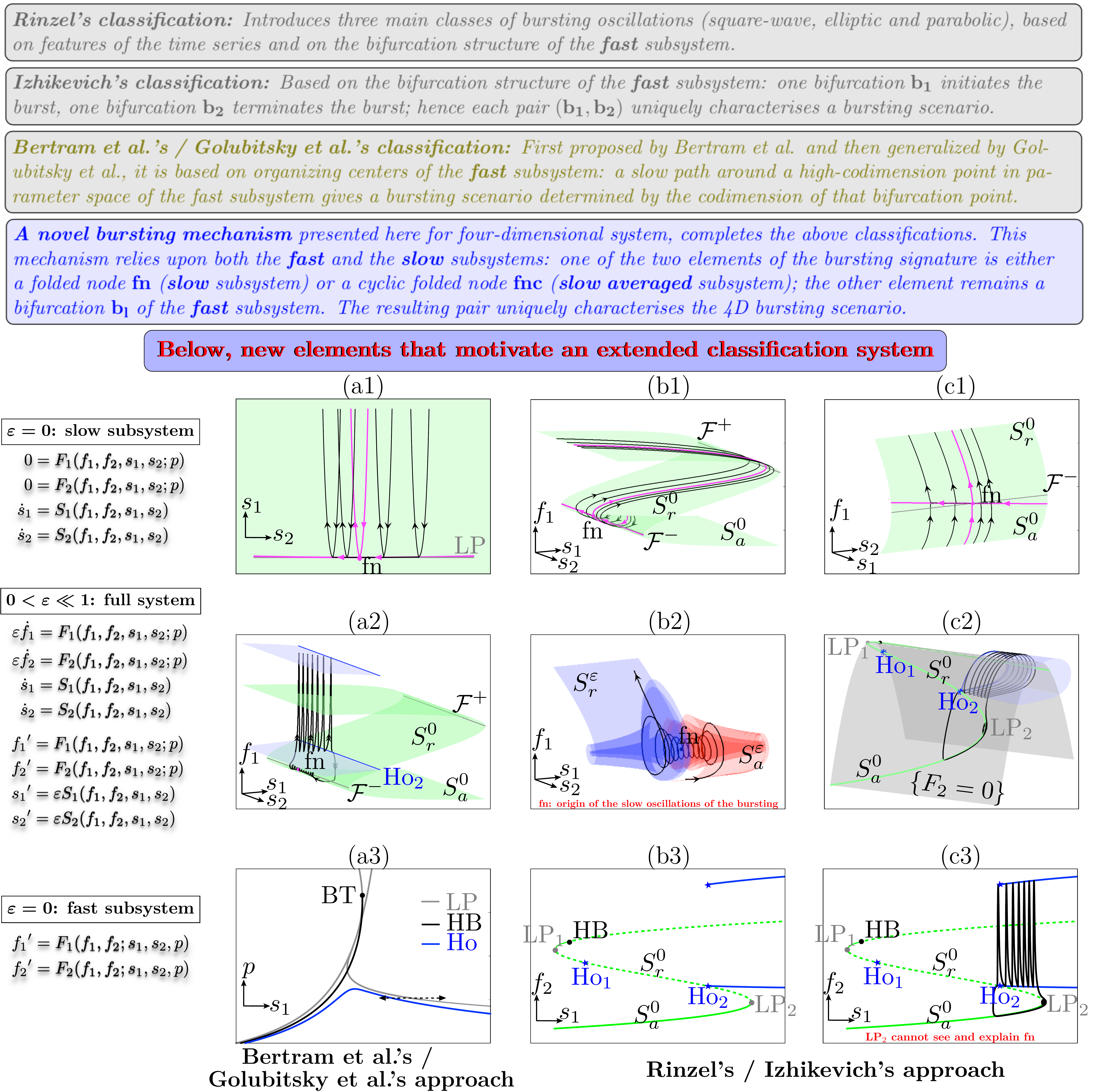}
\caption{(Top part) Rationale behind bursting classifications: that of Rinzel, Izhikevich, Golubitsky et al., respectively, as well as the folded-node bursting approach proposed in this paper. (Main part (ai)-(ci), i=1,2,3) An exemplary folded-node bursting scenario, namely folded-node-homoclinic bursting, is shown in (a2)-(c2) in different 3D phase-space projections; (b2) is a zoomed view of (a2) where the critical manifold $S^0$ has been replaced by attracting $S^a_{\eps}$ and repelling $S^r_{\eps}$ slow manifolds in the vicinity of the folded node (dot). The classical analysis using the approaches of Bertram et al. / Golubitsky et al. (a3), Rinzel and Izhikevich (b3)-(c3) is based upon the fast subsystem's bifurcation structure. The complementary key approach to fully characterize this bursting pattern, presented in panels (a1)-(c1), uses the slow subsystem and its flow which reveals the presence of the folded node \textbf{fn}. \mdsr{Note that the fold curve $\mathcal{F}^-$ of the critical manifold $S^0$ (panel (a2)) is also a curve of saddle-node bifurcation points (labelled LP in (a1)) of the fast subsystem when both slow variables are considered as parameters, the same curve used in the singularity theory approach on panel (a3).} (Side part) Full system's equations as well as the various subsystems, fast and slow, obtained for $\eps=0$.}
\label{fig:sketch}
\end{figure}

It is compelling to ask if there are other bursting oscillations beyond the Rinzel-Izhikevich and Bertram-Golubitsky classification approaches, as summarised in Figure~\ref{fig:sketch} (top panel), and which cannot be explained invoking these state-of-the-art results. If so, could there be an improved classification system that captures a larger class of bursting dynamics beyond the ca. 120 cases captured by the state-of-the-art? Noteworthy, it is reasonable to contemplate that there could be possible extensions for classifying bursting patterns in systems with more than two timescales. In this context, the mathematical analysis would have to deal with nested fast subsystems, which has not yet been achieved (except in very particular cases) and therefore there is substantial work to do in order to extend the state-of-the-art approaches. However, still within the two-timescale framework the question of bursting classification extension remains. This question gains further support since there are electrophysiological recordings of bursting dynamics, which resist the state-of-the-art classification system. A case in point is depicted in Figure~\ref{fig:MMBO_exp}, where the bursting oscillations has two phases, but the quiescent phase has the peculiarity that it appears to periodically rise close to a threshold, however the neuron does not have a transition to the active phase and instead descends back to its baseline activity and only on the second run the active phase emerges. 
These observations suggest that there is an underlying complex mechanism for the quiescent phase of the oscillations and therefore point towards a bursting classification framework that has to also incorporate the analysis of the slow subsystem, which is in stark contrast to state-of-the-art approaches. Further motivating this view is our earlier study, which constructed the first example (to the best of our knowledge) of a slow-fast bursting model whereby the initiation of the bursting oscillations could not be explained by the fast subsystem of the underlying slow-fast model~\cite{desroches13a}. However, therein we did not attempt to derive an improved bursting classification framework. Thus, we herein propose an extension of the state-of-the-art classification system that enforces the importance of considering complex dynamical mechanisms within the slow subsystem as the underlying cause of the initiation or termination of bursting oscillations, which can not be explained by the fast subsystem's analysis solely. To substantiate this novelty we will also show in subsequent sections how to construct a variety of these new cases of bursting oscillations. However, to better guide the reader throughout this manuscript we first replicate the results of our previous work in Figure~\ref{fig:sketch} (panels (a1)-(c3)) and further summarise the key insights of the proposed novel classification framework. We first considered the minimal setting of systems with two-timescales and that possess at least two slow variables. We then constructed a bursting model whose quiescent phase displays small-amplitude near-threshold oscillations. Mathematically, it turns out that these observations are best explained by the so-called \textit{folded-node singularities} defined in the slow subsystem ($\eps=0$) and associated \textit{canard solutions}, which persist for small enough $\eps>0$; see Figure~\ref{fig:sketch} panels (a1)-(c1). Noteworthy, these are not \emph{per se} bifurcations of the fast subsystem but lie on saddle-node bifurcation curves of the fast subsystem (see Figure~\ref{fig:sketch} panel (a1)). 
Such unexpected and non-trivial emergent mathematical objects allow trajectories of the slow subsystem to visit both the attracting ($S^0_a$) and repelling ($S^0_r$) parts of the critical manifold. In the full system (for small $\eps\!>\!0$) the perturbed versions of these manifolds --attracting $S^{\eps}_a$ and repelling $S^{\eps}_r$ slow manifolds-- twist and intersect multiple times (see Figure~\ref{fig:sketch} panels (a2)-(b2)) thereby causing trajectories to non-trivially and robustly oscillate during the quiescent phase. \mdsr{In essence, if the slow dynamics change direction along a fold then this can create a folded node scenario.} The transition from quiescent to active phase is caused by a repulsion of the trajectory away from the unstable sheet of the critical manifold; this phenomenon is mediated by canard solutions due to the presence of the folded node in the slow subsystem. In this particular example, the fast oscillations of the active phase are due to a nearby family of Hopf bifurcations in the fast subsystem (Figure~\ref{fig:sketch} panel (c2)). The return back to quiescence is then caused by a family of homoclinic bifurcations of the fast subsystem. The key insight is that the fast subsystem is blind to what is causing these small-amplitude oscillations during the quiescent phase, and thus it is unable to classify the initiation of these oscillations based upon the bifurcations of fast subsystem only. This point is illustrated by the Rinzel-Izhikevich slow-fast dissection and projection of the trajectory of the full system onto the bifurcation diagram of the fast subsystem (see Figure~\ref{fig:sketch} panels (b3)-(c3)). Note that by employing the  Rinzel-Izhikevich classification system, the bursting dynamics would be explained by two bifurcations of the fast subsystem, namely the fold bifurcation LP$_2$ and the homoclinic bifurcation Ho$_2$. In particular, a fold bifurcation (LP$_2$) does not explain an oscillation. Moreover, a similar argument applies to the Golubitsky approach (see Figure~\ref{fig:sketch} panel (a3)). This panel displays curves of codimension-one bifurcation points of the fast subsystem, which meet at codimension-two e.g. a Bogdanov-Takens BT (within a two-dimensional parameter space). It can then be shown that it is impossible to construct a path for the slow dynamics (within this two-dimensional parameter space), in particular along the homoclinic and saddle-node curves (since these characterise the bursting in the fast subsystem), which could explain folded-node-initiated quiescent phase oscillations. It turns out amongst all possible folded singularities, only folded nodes (and in limiting cases, so-called \textit{folded saddle-nodes}) can generate such robust small-amplitude oscillations in the full system, and this is due to the twisting of nearby attracting and repelling slow manifolds. The key message is that these folded-node singularities only appear in the slow limit of the underlying two-slow-two-fast bursting systems and are invisible in the fast limit. Therefore one must consider both fast and slow subsystems in order to fully characterise the novel bursting scenarios associated with folded-node singularities, which leads us to a novel bursting classification system (see Figure~\ref{fig:sketch} top panel in blue for the new framework). We believe these insights will fuel subsequent developments in higher dimensional multi-scale systems.

As hinted \mdsr{above}, there are emergent slow dynamical mechanisms (captured by slow variables) that are blind to the state-of-the-art classification and thus this suggests the need to extend the classification of bursting oscillations. An example of such emergent slow dynamics is the so-called \textit{folded node} and, in a minimal setting, it emerges due to non-trivial interactions between two slow variables. Indeed, herein we focus on new classes of bursting oscillations modelled via slow-fast systems with (at least) two slow variables and (at least) two fast variables and for which some epochs of the oscillatory time-series is explained by folded-node dynamics. This underlying folded-node signature leads us to name the resulting new classes of bursting models, \textit{folded-node bursters}. Three fundamental cases are envisaged. The first case are bursters characterised by small-amplitude oscillations that occur during the quiescence phase, in which case we will refer to the \textit{classical folded-node} bursting scenario. The second case involves slow-amplitude modulation of the burst, which we will denote as the \textit{cyclic folded-node} bursting scenario. The third case, combines \textit{classical folded-node} and \textit{cyclic folded-node}. These classes of bursting patterns involve both the fast subsystem and the slow subsystem of the model, unlike traditional bursters. A second key aspect of these new classes is the central role played by \textit{canards}, namely, \textit{spike-adding canard} cycles involved in the classical folded-node bursting case, and \textit{torus canards} in the cyclic folded-node bursting case. In the following subsections, we describe in details these two scenarios.

\subsection{Classical folded-node case}
\label{sec:classicalFN}
%
Here we propose several bursting oscillations mediated by a classical folded node and the modelling steps of underlying \mdsr{idealized} models are given. To guide the reader towards a modelling strategy of these systems, we first recall key concepts and mechanisms.

\subsubsection{A necessary preliminary step: spike-adding canard explosion}
\label{sec:spikeadd}
First, recall that \textit{canards} are non-trivial trajectories that emerge due timescale separation and unexpectedly, these trajectories contain segments that follow both an attracting and a repelling slow manifold. This phenomenon has been thoroughly studied in planar systems (i.e. with 1 slow variable and 1 fast variable)~\cite{benoit81,dumortier96,eckhaus83,krupa01,mishchenko94}, as well as, in 3D systems with two slow variables~\cite{benoit90,desroches12,wechselberger05}. In applications, canards can be associated to complex (bio)physical mechanisms, for example in neuroscience it provides the best approximation to the excitability threshold in certain single-neuron models. This observation was first made by Izhikevich~\cite{izhikevich07}, who showed that canards organise the transition to the spiking regime of type II neurons. This was later analysed in more details in~\cite{demaesschalck15,desroches13b,wechselberger13}. Yet another important mechanism is the so-called \textit{spike-adding canard explosion}. This canard phenomenon arises in bursting oscillations and can be described as a sequence of canard explosions which organise the transition from subthreshold oscillations to bursting solutions with more and more spikes per bursts. This phenomenon was first described and analysed (in the case of chaotic dynamics regime) in~\cite{terman91} in the context of square-wave bursting. This was revisited more recently in~\cite{guckenheimer09} from the computational standpoint of saddle-type slow manifolds and further described in~\cite{nowacki12} in a modeling context to explain transient spikes. These analyses were later refined (from a canard standpoint) in~\cite{desroches13a} and the canard-mediated spike-adding dynamics was fully analysed in~\cite{desroches16} in the context of parabolic bursters (with two slow variables), revealing the central role of folded-saddle canards. Noteworthy, bursting oscillations that possess spike-adding mechanism is a limiting (border line) case that already hints to the importance of possibly including the analysis of the slow flow in a bursting classification framework. That is, spike-adding requires a turning point of the slow flow, a so-called \textit{canard point}, whereby each new added spike (within the bursting phase) is born via a slow (delayed) passage through this turning point. Crucially, the fast subsystem is blind to the underlying canard trajectories occurring near the turning point (well-defined as such only in the slow flow) and instead only sees a fold bifurcation. Therefore the state-of-the-art bursting classification systems does not capture this aspect. Nevertheless, we refrain from declaring this phenomenon as a new bursting mechanism because a spike-adding canard explosion gives rise to canard cycles that exist only within exponentially thin parameter regions. Hence, the robust dynamics is the fold-initiated bursting dynamics, and the fast subsystem analysis still prevail in order to classify it. In contrast, if we consider a fold-initiated bursting scenario undergoing spike-adding canard explosion and if we further add a slow dynamics for the parameter that displays the spike-adding canard explosion (i.e. a second slow variable in the extended model) then we obtain a folded-node bursting system. This has a similar effect to the case in classical (van der Pol type) systems where the canard phenomenon becomes robust if one adds a second slow variable, which has the effect of creating a folded singularity in the resulting two-dimensional slow flow and allows for multiple canard trajectories to exist. The idea here is similar, but with two fast variables, allowing for bursting dynamics in conjunction with folded-node dynamics. A first example of this scenario was termed \textit{mixed-mode bursting oscillations} in~\cite{desroches13a} but we prefer to denote it more generally folded-node bursting. Indeed, folded-node bursting is a new form of bursting pattern with two slow variables where the silent phase contains small-amplitude (subthreshold) oscillations due to the presence of a folded node in the slow subsystem. This folded node is responsible for the presence of a funnel region in the full system and trajectories entering this funnel make a number of rotations (which can be controlled by adjusting parameters) before they leave it and start to burst. Hence, the passage through the folded-node funnel organises the transition from quiescence to burst and it can only be understood by suitably analysing the slow subsystem. We subsequently describe a strategy for constructing folded-node bursting systems.
\subsubsection{Construction of minimal folded-node bursting systems}
As a staring point, we consider the prototypical fold-initiated burster of Hindmarsh-Rose type. By this we mean a three-dimensional slow-fast system with two fast variables and one slow and a cubic-shaped family of equilibria in the fast subsystem, namely the critical manifold $S^0$. We can write the following set of differential equations (using the fast time $\tau$) to describe the dynamics of such a system
\begin{equation}\label{eq:protoburster}
\begin{split}
x' &= y - f(x) + az,\\
y' &= g(x)-y,\\
z' &= \eps(\alpha x + \gamma\beta - \delta z),
\end{split}
\end{equation}
where $f$ is a cubic polynomial function, $g$ is (at least) quadratic; moreover, $0\!<\eps\!\ll1$ is a small parameter and $(a,\alpha,\beta,\gamma,\delta)$ are potential bifurcation parameters; why we use a product of two parameters in the $z$ equation will become clear in the next section. As we shall see in the example section to follow, one can also obtain all fold-initiated scenarios by using an unfolding of a codimension-3 degenerate Bogdanov-Takens (BT) bifurcation; see~\cite{dumortier91} for details.

\begin{figure}[!t]
\centering
\includegraphics{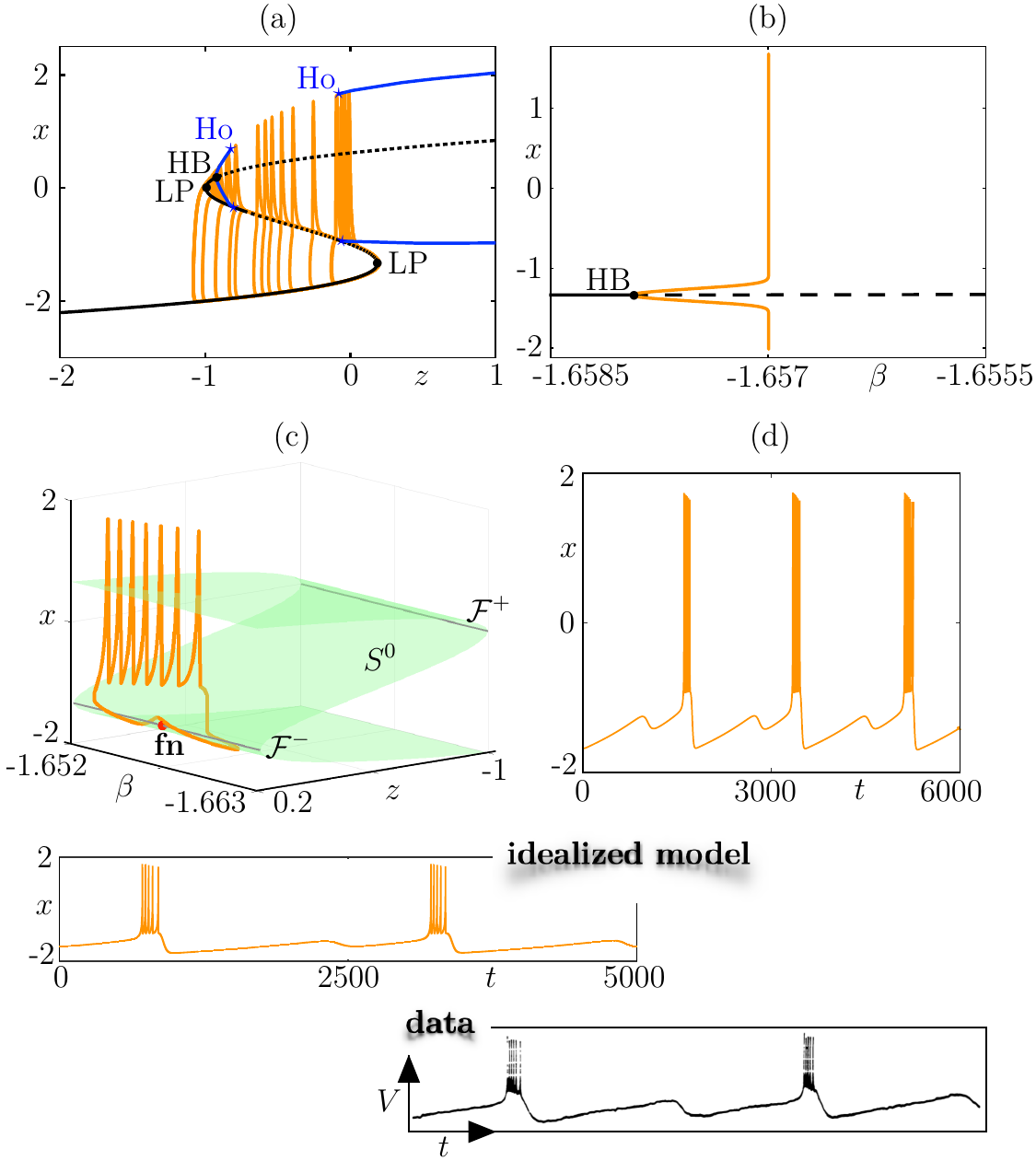}
\caption{Folded-node/Homoclinic bursting. We take $f(x)=x^3-3x^2$ and \bluebis{$G(x,y,z)=g(x)-y$ with $g(x)=1-5x^2$. The parameter values are: $a=1$, $c=1$, $\alpha=0.3$, $\gamma=1$, $\delta=1.2$, $\eps=0.002$, $\mu=0.033$, $\gamma_y=0.0005$ and $\gamma_{\beta}=-0.008$}. Panels (a-b) show the spike-adding transition in system~\eqref{eq:fnburster}: (a) in the $(z,x)$ plane; (b) associated bifurcation diagram with respect to parameter $\beta$. Panels (c-d) show a folded-node/homoclinic bursting orbit: (c) in the $(\beta,z,x)$ space; (d) $x$-time series. The bottom panels show a comparison between this folded-node bursting orbit from~\eqref{eq:fnburster} and experimental data from~\cite{roy84}.}
\label{fig:fnhom}
\end{figure}

\begin{figure}[!t]
\centering
\includegraphics{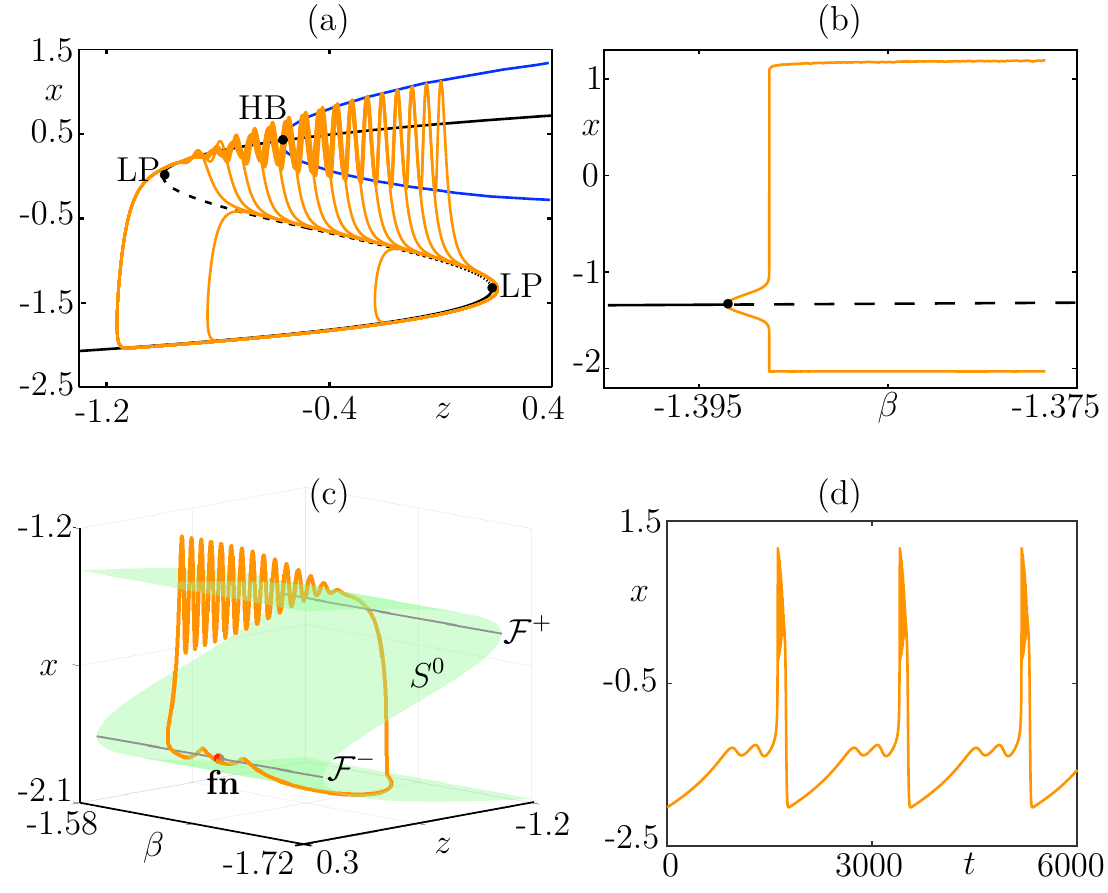}
\caption{Folded-node/Hopf bursting. We take $f(x)=x^3-3x^2$ and \bluebis{$G(x,y,z)=g(x)-y$ with $g(x)=1-5x^2$. The parameter values are: $a=1$, $c=2$, $\alpha=0.3$, $\gamma=1$, $\delta=1$, $\eps=0.004$, $\mu=0.0104$, $\gamma_y=0.0003$ and $\gamma_{\beta}=-0.05$}. Panels (a-b) show the spike-adding transition in system~\eqref{eq:fnburster}: (a) in the $(z,x)$ plane; (b) associated bifurcation diagram with respect to parameter $\beta$. Panels (c-d) show a folded-node/Hopf bursting orbit: (c) in the $(\beta,z,x)$ space; (d) $x$-time series.}
\label{fig:fnhopf}
\end{figure}

\begin{figure}[!t]
\centering
\includegraphics{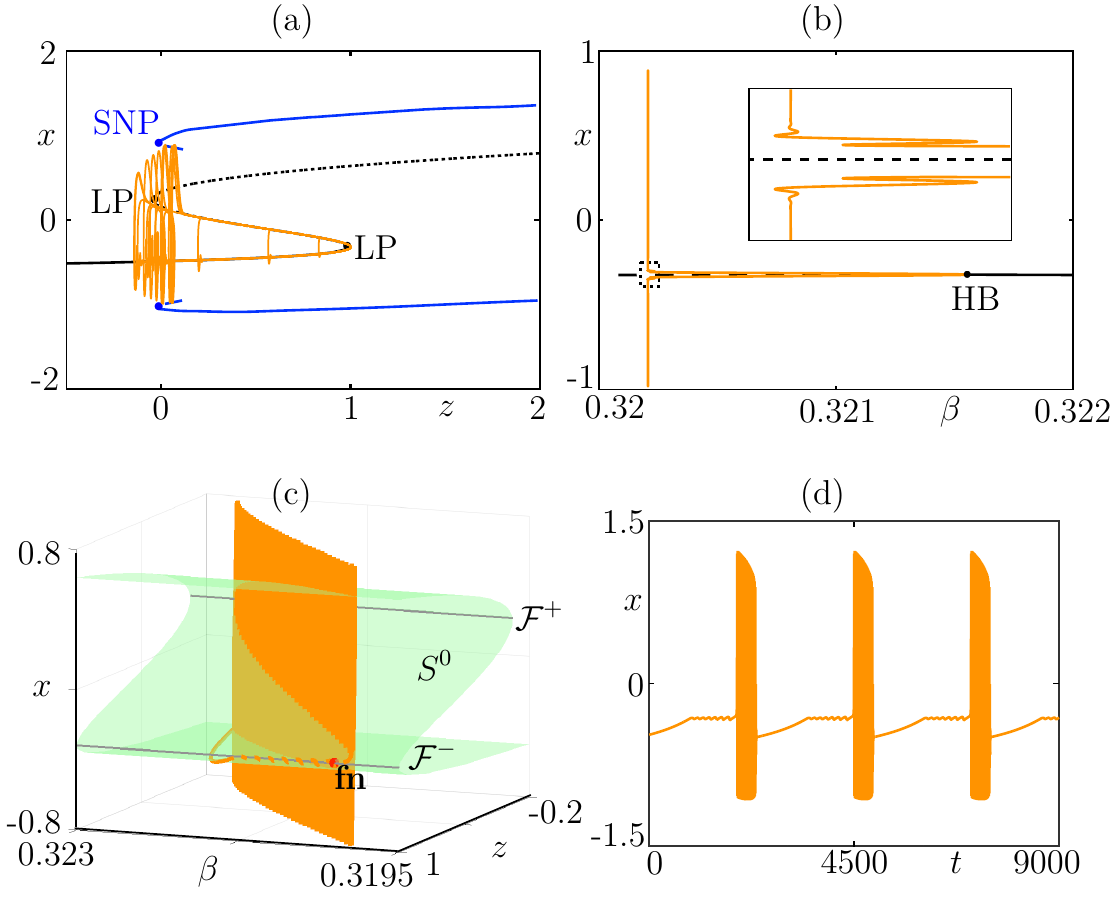}
\caption{Folded-node/Fold of cycles bursting. We take $f(x)=0$ and \bluebis{$G(x,y,z)=-x^3+A_1(z)x+A_2(z)-y(A_3(z)-x+x^2)$}, where \bluebis{$A_1(z)=0.1201z+0.1871$, $A_2(z)=0.0906z-0.0251$, $A_3(z)=0.105z-0.3526$}. \bluebis{The parameter values are: $a=0$, $c=1$, $\alpha=0$, $\gamma=-1$, $\delta=1$, $\eps=0.01$, $\mu=-0.00012$, $\gamma_y=-0.003$, $\gamma_{\beta}=0.0001$}. Panels (a-b) show the spike-adding transition in system~\eqref{eq:fnburster}: (a) in the $(z,x)$ plane; (b) associated bifurcation diagram with respect to parameter $\beta$. Panels (c-d) show a folded-node/fold of cycles bursting orbit: (c) in the $(\beta,z,x)$ space; (d) $x$-time series.}
\label{fig:fnsnp}
\end{figure}

A few assumptions are required in order for the system~\eqref{eq:protoburster} to display fold-initiated bursting. First of all, we assume that $f$ and $g$ are adequately chosen so that the fast subsystem has a cubic-shaped family of equilibria that depends on $z$ as a parameter (for the fast subsystem). This family can be written as a cubic function: 
$$z=\left(f(x)-g(x)\right)/a.$$ 
Therefore, the corresponding bifurcation diagram (of the fast subsystem) in $z$ is S-shaped and will have fold points. 
The critical manifold is then given by 
\begin{eqnarray}\label{eq:critman}
S^0:=\left\{(x,y,z)\in\mathbb{R}^3\;\big/\;\;y=g(x)\;,\;z=\left(f(x)-y\right)/a\right\}.
\end{eqnarray}
We also require bistability in the fast subsystem between equilibria and limit cycles, in an interval of $z$-values. One bound of this interval correspond to a fold bifurcation and, geometrically, to one fold point of the cubic family of equilibria. The other boundary of the region of bistability of the fast subsystem will be a bifurcation of limit cycles and we shall consider three main cases, namely, saddle-homoclinic bifurcation (see Fig.~\ref{fig:fnhom}), Hopf bifurcation (see Fig.~\ref{fig:fnhopf}) and fold bifurcation of cycles (see Fig.~\ref{fig:fnsnp}), but the list is not exhaustive. Now, considering the linear slow dynamics of system~\eqref{eq:protoburster} for the slow variable $z$, we assume that a variation of one of the two parameters $\alpha$ and $\beta$ in the full system induces the linear $z$-nullsurface to cut through the fold point of the critical manifold $S^0$ for a certain value of this parameter. One can show that this creates a Hopf bifurcation in the full system which induces limit cycles to appear. Provided this transversal cut of the $z$-nullsurface with the critical manifold takes place, then a spike-adding canard explosion will emerge, whereby bursting solutions appear from subthreshold (spikeless) periodic solutions along branch of limit cycles undergoing multiple canard explosions; see~\cite{desroches13a} for an example of this phenomenon in the context of square-wave bursting. As explained in the previous section, one salient feature of the spike-adding canard explosion is the presence of a turning point (a canard point) in the slow flow of system~\eqref{eq:protoburster}. To compute the slow-flow, we first rescale time in~\eqref{eq:protoburster} by a factor $\eps$. That is, we rescale the fast time $\tau$ (with $x'=dx/d\tau$) into the slow time $t$ defined by $t=\eps\tau$. This bring the system to the slow-time parametrisation 
\begin{equation}\label{eq:protobursterslowtime}
\begin{split}
\eps\dot{x} &= y - f(x) + z,\\
\eps\dot{y} &= g(x)-y,\\
~~\dot{z} &= (\alpha x + \gamma\beta - \delta z),
\end{split}
\end{equation}
whose $\eps=0$ limit corresponds to the slow subsystem. The slow subsystem is a differential-algebraic equation (DAE), where the dynamics of $z$ is explicitly preserved while $x$ and $y$ are slaved to $z$ by the algebraic constraints that corresponds to the equation of the (here one-dimensional) critical manifold $S^0$. The dynamics of $x$ and $y$ can be revealed by differentiating the algebraic constraint with respect to the slow time, which gives after rearranging the following one-dimensional dynamical system defined on $S^0$
\begin{equation}\label{eq:protobursterslowflow}
\begin{split}
\dot{x} &= a\frac{\alpha x + \gamma\beta - \delta z}{f'(x)-g'(x)}.
\end{split}
\end{equation}
%
As is \red{typical} in slow-fast systems with folded critical manifolds, note that the denominator of the right-hand side of~\eqref{eq:protobursterslowflow} vanishes at fold points of $S^0$, which makes \mdsr{generically} the dynamics of $x$ explode \mdsr{and the corresponding fold point is referred to as a \textit{jump point}. However, if} the numerator has a zero of the same order \mdsr{as the denominator, then there can be a cancellation and the dynamics of $x$ does not explode; in this case, the fold point is referred to as a \textit{canard point} or a \textit{turning point}. The condition for a canard point to occur \red{in this system is}} \mdsr{given} by the following condition: 
\begin{equation}\label{eq:protoburstercanardpoint}
\begin{split}
z_f &= (\alpha x_f + \gamma\beta)/\delta,
\end{split}
\end{equation}
where $(x_f,z_f)$ is a fold point of $S^0$. This indeed gives a transversal crossing of the slow nullsurface with the critical manifold at one of its fold points. Even though~\eqref{eq:protoburstercanardpoint} depends on several parameters, it is a codimension-1 condition, therefore by fixing all parameters but one, then the condition can be satisfied by adjusting the last parameter. We arbitrarily choose to vary $\beta$, which will become a second slow variable in the full 4D folded-node bursting system that we will construct below. Therefore, the spike-adding transitions leading to bursting in system~\eqref{eq:protoburster} are obtained as the result of the slow nullsurface moving though one fold point of the critical manifold upon variation of $\beta$. The same dynamics would be obtained by varying a parameter affecting the critical manifold while maintaining the slow nullsurface fixed, in particular if we were to append an additive parameter $I$ to the $x$ equation of the system. This would mimick the effect of an applied (external) current in a neuron-type model such as the Hindmarsh-Rose model~\cite{hindmarsh84} or the Morris-Lecar model~\cite{morris81,terman91}. However, from the pure dynamical viewpoint, varying a parameter in the slow equation results in the same effect and this is the scenario that we chose in order to construct fold-initiated spike-adding transitions in the original 3D burster and folded-node bursting in the extended 4D model.

Starting from a fold-initiated bursting scenario with spike-adding canard explosion (controlled via a static variation of parameter $\beta$) then a folded-node bursting is obtained by prescribing the dynamics on $\beta$ by a slow differential equation. That is, we consider the following extended bursting system
\bluebis{
\begin{equation}\label{eq:fnburster}
\begin{aligned}
x' &= (y - f(x) + az)/c,\\
y' &= G(x,y,z),\\
z' &= \eps(\alpha z + \gamma\beta - \delta x),\\
\beta'  &= \eps\big(\mu-\gamma_y(y-y_{\mathrm{fold}})^2-\gamma_{\beta}(\beta-\beta_{\mathrm{fold}})^2\big).
\end{aligned}
\end{equation}
}
\bluebis{Note that we will consider prototype systems where either $G$ is directly given as a graph over $x$, described as, $G(x,y,z)=g(x)-y$ (i.e. Folded-node/homoclinic and folded-node/Hopf cases), or the level set $\{G(x,y,z)=0\}$ is a graph over $(x,z)$, expressed as, $\{y=g(x,z)\}$ (i.e. Folded-node/fold of cycles case). We claim that all folded-node initiated bursting scenarios can be obtained in either of these two ways. In the latter case, our minimal model is inspired by the codimension-3 degenerate Bogdanov-Takens unfolding introduced in~\cite{dumortier91} and further applied in the context of bursting in~\cite{saggio17}. In practice (for simulation purposes), $\mu$, $\gamma_y$ and $\gamma_z$ will be taken $O(\eps)$.} Therefore, system~\eqref{eq:fnburster} is effectively a three-timescale dynamical systems with dynamics evolving on $O(1)$, $O(\eps)$ and $O(\eps^2)$ timescales. For convenience and to ease the folded-node analysis, we will keep the equations written as in~\eqref{eq:fnburster} with only $\eps$ has an apparent timescale separation parameter. System~\eqref{eq:fnburster} is parametrised by the fast time $\tau$, meaning that the $\eps=0$ limit of that parametrisation of the system gives the fast subsystem where the slow variables $z$ and $I$ are frozen and become parameters. Introducing the slow time $t=\eps\tau$ brings the system into the different parametrisation 
\bluebis{
\begin{equation}\label{eq:fnburstertau}
\begin{aligned}
\eps\dot{x} &= (y - f(x) + az)/c,\\
\eps\dot{y} &= G(x,y,z),\\
\dot{z} &= \alpha z + \gamma\beta - \delta x,\\
\dot{\beta}  &= \mu-\gamma_y(y-y_{\mathrm{fold}})^2-\gamma_{\beta}(\beta-\beta_{\mathrm{fold}})^2,
\end{aligned}
\end{equation}
}
whose $\eps\!=\!0$ limit corresponds to the slow subsystem. We will show that, all other parameters being fixed, the slow subsystem of~\eqref{eq:fnburstertau} possesses a folded-node singularity, which creates transient subthreshold oscillations that initiate the burst when $0<\eps\ll1$, regardless of the values of other parameters. However, simulations require that $\mu$\, $\gamma_y$ and $\gamma_I$ be $O(\eps)$ in order for these small subthreshold oscillations to be recurrent, hence entering into a robust periodic bursting attractor which we name folded-node bursting. We provide numerical evidence of this point, based on the strength of the global return mechanism, even though we do not provide a rigorous proof of it.

\bluebis{Applying the same strategy as in the three-dimensional (bursting) case, and projecting onto the $(x,\beta)$-plane (the dimension of the slow flow corresponds to the number of slow variables), we obtain the following system for the reduced system (or slow subsystem)
\begin{equation}\label{eq:fnbursterslowflow}
\begin{aligned}
\dot{x} &= \frac{(g_z(x,z)+a)(\alpha z + \gamma\beta - \delta x)}{f'(x)-g_x(x,z)},\\
\dot{\beta}  &= \mu-\gamma_y\left(g(x,z)-y_{\mathrm{fold}}\right)^2-\gamma_{\beta}(\beta-\beta_{\mathrm{fold}})^2,
\end{aligned}
\end{equation}
after substituting for $g(x,z)$ for $y$ from the critical manifold condition. Note that in two of the three examples that we will consider, $g(x,z)=g(x)$ depends only on $x$ and hence $g_z(x,z)=0$. The critical manifold of system~\eqref{eq:fnburster} is not normally hyperbolic everywhere and, hence, the system possesses a (1D here) fold set defined by $$\mathcal{F}:=\{(x,y,z)\in S^0;\;f'(x)=g_x(x,z)\}.$$
This implies that the slow flow~\eqref{eq:fnbursterslowflow} of system~\eqref{eq:fnburster} is not defined along $\mathcal{F}$. The slow flow can be extended along $\mathcal{F}$ by preforming an $x$-dependent time rescaling which amounts to multiply the right-hand side of~\eqref{eq:fnburster} by a factor $f'(x)-g_x(x,z)$, hence yielding the so-called \textit{desingularised reduced system (DRS)}
\begin{equation}\label{eq:fnbursterDRS}
\begin{aligned}
\dot{x} &= \left(g_z(x,z)+a\right)(\alpha z + \gamma\beta - \delta x),\\
\dot{\beta}  &= \left(f'(x)-g_x(x,z)\right)\left(\mu-\gamma_y\left(g(x,z)-y_{\mathrm{fold}}\right)^2-\gamma_{\beta}(\beta-\beta_{\mathrm{fold}})^2\right),
\end{aligned}
\end{equation}
with $z=z(x)$ defined by $S^0$, that is, $g(x,z)-f(x)+\alpha z=0$. In all cases we shall consider (including the general codimension-3 unfolding of a degenerate BT bifurcation from~\cite{dumortier91}), $z$ can be written as a function of $x$ on $S^0$. 
As a consequence of this $x$-dependent time rescaling, the DRS~\eqref{eq:fnbursterDRS} is regular everywhere in $\mathbb{R}^2$ including on $\mathcal{F}$, along which it has the possibility for equilibria simply by appearance of the factor $f'(x)-g_x(x,z)$ in the $\beta$-equation. The equilibrium condition is then that $\dot{x}=0$ in~\eqref{eq:fnbursterDRS} together with $f'(x)-g_x(x,z)$, which conveys the idea already seen in the 3D (bursting) case. That is, a singularity of the reduced system at a point on $\mathcal{F}$ can be resolved if and only if the numerator of the right-hand side of $\dot{x}$ in that system vanishes at this point and the zeros of the two algebraic expressions to be of the same order. Such points are called \textit{folded singularities} (or \textit{folded equilibria}) and they are the equivalent of canard points in the cases with (at least) two slow variables. Folded equilibria are equilibria of the DRS~\eqref{eq:fnbursterDRS} and, according to their topological type as equilibria of the DRS, one can generically define \textit{folded nodes}, \textit{folded saddles} and \textit{folded foci}. However they are not equilibria of the reduced system~\eqref{eq:fnbursterslowflow} due to the singular time rescaling performed to pass from one to the other. Indeed, this time rescaling is chosen so that trajectories of the DRS have reversed orientation on the repelling sheet of $S^0$ compared to trajectories of the reduced system (both have the same orientation along the attracting sheet). Hence, in the case of folded nodes and folded saddles, trajectories starting on the attracting sheet of $S^0$ may cross the folded singularity in finite time and with finite speed.\newline
The Jacobian matrix of~\eqref{eq:fnbursterDRS} evaluated at a folded equilibrium has the form
\begin{equation}\label{eq:JacDRS}
\begin{aligned}
\mathrm{J}=\begin{pmatrix}
(-\delta+\alpha z'(x))(g_z(x,z)+a) & \gamma(g_z(x,z)+a)\\
K_2 & 0
\end{pmatrix},
\end{aligned}
\end{equation}
where $$K_2=(f''(x)-\partial_xg_x(x,z))\left(\mu-\gamma_y\left(g(x,z)-y_{\mathrm{fold}}\right)^2-\gamma_{\beta}(\beta-\beta_{\mathrm{fold}})^2\right).$$ From~\eqref{eq:JacDRS}, one can easily write down conditions that enable the emergence of a folded-node singularity ($\mathrm{tr}(J)<0$, $\det{J}>0$, $\mathrm{tr}(J)^2-4\det{J}>0$) or a folded-saddle singularity ($\det{J}<0$) in the reduced system. As we will explain below, even though only the folded-node case gives rise to robust bursting patterns, the folded-saddle case is still interesting in the study of 4D bursters with two slow variables. One also can easily verify that our minimal example systems all give rise to a folded-node case. Indeed, in the folded-node/homoclinic (Figure~\ref{fig:fnhom}) and folded-node/Hopf (Figure~\ref{fig:fnhopf}) bursting cases, system~\eqref{eq:fnburster} has the form
\begin{equation}\label{eq:fnhomburster}
\begin{aligned}
x' &= (y - x^3 + 3x^2 + z)/c,\\
y' &= 1 - 5x^2 - y,\\
z' &= \eps(\alpha z + \gamma\beta - \delta x),\\
\beta'  &= \eps\left(\mu-\gamma_y(y-y_{\mathrm{fold}})^2-\gamma_{\beta}(\beta-\beta_{\mathrm{fold}})^2\right),
\end{aligned}
\end{equation}
which hence gives the following DRS's Jacobian matrix
\begin{equation}\label{eq:JachomDRS}
\begin{aligned}
\mathrm{J_{1,2}}=\begin{pmatrix}
-\delta & \gamma\\
K_2 & 0
\end{pmatrix},
\end{aligned}
\end{equation}
with: $K_2=(-6x_{\mathrm{fs}}-4)\left(\mu-\gamma_y\left(1-5x_{\mathrm{fs}}^2-y_{\mathrm{fold}}\right)^2-\gamma_{\beta}(\beta-\beta_{\mathrm{fold}})^2\right)$, and $x_{\mathrm{fs}}=-4/3$. Given the chosen parameter values corresponding to Figures.~\ref{fig:fnhom} and~\ref{fig:fnhopf}, then we immediately conclude that we have indeed a folded node. Likewise, in the folded-node/fold of cycles case illustrated in Figure~\ref{fig:fnsnp}, the slow-fast system corresponding to~\eqref{eq:fnburster} is 
\begin{equation}\label{eq:fnsnpburster}
\begin{aligned}
x' &= y,\\
y' &= -x^3+A_1(z)x+A_2(z)-y(A_3(z)-x+x^2),\\
z' &= \eps(\alpha z + \gamma\beta - \delta x),\\
\beta'  &= \eps\left(\mu-\gamma_y(y-y_{\mathrm{fold}})^2-\gamma_{\beta}(\beta-\beta_{\mathrm{fold}})^2\right),
\end{aligned}
\end{equation}
where $A_i=a_i z+b_i$ ($i=1,2,3$) are linear functions of $z$. Therefore, we obtain the associated DRS's Jacobian matrix
\begin{equation}\label{eq:JacsnpDRS}
\begin{aligned}
\mathrm{J_{1,2}}=\begin{pmatrix}
\left(-\delta+\alpha\frac{3x_{\mathrm{fs}}^2-b_1}{a_1+a_2}\right)(a_1x_{\mathrm{fs}}+a_2) & \gamma(a_1x_{\mathrm{fs}}+a_2)\\
K_2 & 0
\end{pmatrix},
\end{aligned}
\end{equation}
with: $K_2=(6x_{\mathrm{fs}}-a_1)\left(\mu-\gamma_yy_{\mathrm{fold}}^2-\gamma_{\beta}(\beta-\beta_{\mathrm{fold}})^2\right)$ and $x_{\mathrm{fs}}$ solution to $$-3x_{\mathrm{fs}}^2+a\frac{x_{\mathrm{fs}}^3-b_1x_{\mathrm{fs}}-b_2}{a_1+a_2}+b_1=0.$$ Substituting the parameter values for their chosen numerical value mentioned in the caption of Figure~\ref{fig:fnsnp} allows to conclude that we are indeed dealing with a folded node.}\newline
One can obtain the general DRS~\eqref{eq:fnbursterDRS} by applying implicit differentiation to one algebraic equation only (the right-hand side of the $\dot{x}$ equation in the original system) and substituting $g(x,z)$ for $y$ (coming from the second algebraic equation). This gives the same result as the DRS obtained from both algebraic constraint together. Indeed, in all generality, applying implicit differentiation to the two algebraic equations of the slow subsystem gives 
\begin{equation}\label{eq:RS2fast}
\begin{aligned}
\begin{pmatrix}
-f'(x)        & 1 \\
-g_x(x,z) & 1
\end{pmatrix}
\begin{pmatrix}
\dot{x} \\
\dot{y}
\end{pmatrix}
&=
\begin{pmatrix}
-a \\
g_z(x,z)
\end{pmatrix}
(\alpha z + \gamma\beta - \delta x)\\
\dot{z} &= \alpha z + \gamma\beta - \delta x,\\
\dot{\beta}  &= \mu-\gamma_y(y-y_{\mathrm{fold}})^2-\gamma_{\beta}(\beta-
\beta_{\mathrm{fold}})^2,
\end{aligned}
\end{equation}
which by Kramer's rule is equivalent, after posing $$\mathrm{J}=
\begin{pmatrix}
-f'(x)        & 1 \\
-g_x(x,z) & 1
\end{pmatrix}
,$$ (Jacobian matrix of the original vector field with respect to the fast variables at $\eps=0$) to 
\begin{equation}\label{eq:RS2fast2}
\begin{aligned}
\det(\mathrm{J})
\begin{pmatrix}
\dot{x} \\
\dot{y}
\end{pmatrix}
&= \mathrm{Adj(J)}\begin{pmatrix}
a \\
-g_z(x,z)
\end{pmatrix}
(\alpha z + \gamma\beta - \delta x)\\
\dot{z} &= \alpha z + \gamma\beta - \delta x,\\
\dot{\beta}  &= \mu-\gamma_y(y-y_{\mathrm{fold}})^2-\gamma_{\beta}(\beta-
\beta_{\mathrm{fold}})^2,
\end{aligned}
\end{equation}
where $\det(\mathrm{J})=g_x(x,z)-f'(x)$ and 
$$\mathrm{Adj}(\mathrm{J})=
\begin{pmatrix}
1            & -1\\
g_x(x,z) & -f'(x)
\end{pmatrix},$$ denote the determinant and the adjugate matrix of $\mathrm{J}$, respectively. The previous system is singular when $\det(\mathrm{J})$ vanishes, which happens on the fold set. It can be desingularized by rescaling time by a factor $\det(\mathrm{J})$, which brings the DRS in its most general form, namely
\begin{equation}\label{eq:DRS2fast}
\begin{aligned}
\begin{pmatrix}
\dot{x} \\
\dot{y}
\end{pmatrix}
&= \mathrm{Adj(J)}\begin{pmatrix}
a \\
-g_z(x,z)
\end{pmatrix}
(\alpha z + \gamma\beta - \delta x)\\
\dot{z} &= \det(\mathrm{J})(\alpha z + \gamma\beta - \delta x)\\
\dot{\beta}  &= \det(\mathrm{J})\big(\mu-\gamma_y(y-y_{\mathrm{fold}})^2-\gamma_{\beta}(\beta-
\beta_{\mathrm{fold}})^2\big).
\end{aligned}
\end{equation}
After being projected onto the $(x,\beta)$-space, system~\eqref{eq:DRS2fast} then takes the form
\begin{equation}\label{eq:DRS2fast2}
\begin{aligned}
\dot{x} &=(a+g_z(x,z))(\alpha z + \gamma\beta - \delta x)\\
\dot{\beta} &= (f'(x)-g_x(x,z))\big(\mu-\gamma_y(y-y_{\mathrm{fold}})^2-\gamma_{\beta}(\beta-
\beta_{\mathrm{fold}})^2\big),
\end{aligned}
\end{equation}
which indeed agrees with~\eqref{eq:fnbursterDRS}.

With the above analysis, we can construct in principle any folded-node burster of our liking. We showcase three examples. First, a folded-node homoclinic burster is presented in Fig.~\ref{fig:fnhom} and we also show in the bottom panels a comparison with data from~\cite{roy84} (also displayed in Fig.~\ref{fig:MMBO_exp}). Note that our \mdsr{idealized} model was not initially designed to explain these data, yet the time profiles match remarkably well. The strong similarity between our \mdsr{idealized} model and these data suggest that folded-node bursting constructions could potentially inform the design of biophysical models. Then, a folded-node Hopf burster is presented in Fig.~\ref{fig:fnhopf}. Last, a folded-node fold-cycle burster is shown in Fig.~\ref{fig:fnsnp}. In each of these three figures, we show in panel (a) the classical Rinzel dissection between the bifurcation diagram of the 2D fast subsystem of the underlying 3D fold-initiated bursting model, and several limit cycles of this 3D bursting model undergoing a spike-adding canard explosion. In panel (b), we show the bifurcation diagram of the 3D bursting system upon variation of the parameter which will become the second slow variable of the folded-node burster, and this diagram displays spike-adding canard explosion. Note that a full analysis of these spike-adding canard explosion scenarios is beyond the scope of the present work as it is by-and-large an open research question. In panels (c), we show a folded-node bursting cycle in a 3D phase-space projection together with the critical manifold $S^0$. In panel (d), we plot the time series of this cycle for the fast variable $x$. More examples of folded-node bursting scenarios can be constructed by following the procedure highlighted above and by choosing a different bifurcation of the fast subsystem ending the burst.

Finally, we quickly reflect on why folded-saddle bursting is not robust. The folded-saddle case is simply a different parameter regime in the slow subsystem, however the resulting dynamics is substantially different than that generated by a folded node. In neuron models with (at least) two slow variables, folded saddles and their associated canard solutions play the role of firing threshold. In particular, in the context of bursting system, they have recently been shown to organise the spike-adding transition in parabolic bursters~\cite{desroches16,desroches18}. Counterintuitively, small-amplitude oscillations can also emerge in the vicinity of a folded saddle; see~\cite{mitry17} for a rigorous analysis of this phenomenon and also~\cite{desroches16b,desroches18} for further related work. However, there is no funnel near a folded saddle and the canard dynamics is hence not robust, which applies no matter how many fast variables the system possesses, so in particular in the context of bursting. This is why, in systems with (at least) two fast and two slow variables, only the folded-node case gives rise to a new class of bursting oscillations.

\subsubsection{Existence of folded-node bursting solutions}
Mixed-mode oscillations have been the subject of intense research in the past few decades, in particular (more recently) in the context of multiple-timescale systems; see~\cite{desroches12}. However, there are very few results guaranteeing the existence of MMOs. One example is due to Br{\o}ns, Krupa and Wechselberger in~\cite{brons06}, where they prove the existence of ``simple'' patterns of MMOs ---with one large-amplitude oscillations and $s$ small-amplitude oscillations, hence denoted a $1^s$ MMO--- using a perturbation argument (in $\eps$) from a singular orbit that they construct using both slow and fast limits of the original system. Their strategy can be adapted to prove the existence of a folded-node bursting periodic orbit by $\eps$-perturbation of a singular orbit. The segment of that singular orbit constructed using the slow flow is basically the same as in the MMO case, that is, it lies in the singular funnel or on the strong canard of the folded node. However, in the case of folded-node bursting the fast subsystem is of dimension 2 and displays limit cycles in the fast part of the cycle. Therefore one needs to concatenate the slow-flow segment with a segment along the \textit{averaged slow flow} of the system; see e.g.~\cite{roberts15} for details.  The averaged slow flow is defined using the fast subsystem in its oscillatory regime and averaging out the fast variables $x$ and $y$ along one cycle to define a slow motion for the slow variables $z$ and $\beta$; its equations are hence given by:
\begin{equation}\label{eq:avslowflow}
\begin{aligned}
x' &= (y - f(x) + az)/c,\\
y' &= G(x,y,z),\\
z' &= \alpha z + \gamma\beta - \delta \langle x \rangle,\\
\beta'  &= \mu-\frac{\gamma_y}{T_z}\int_0^{T(z)}(y(s)-y_{\mathrm{fold}})^2 -\gamma_{\beta}(\beta(s)-\beta_{\mathrm{fold}})^2 ds.
\end{aligned}
\end{equation}

\begin{figure}[!t]
\centering
\includegraphics{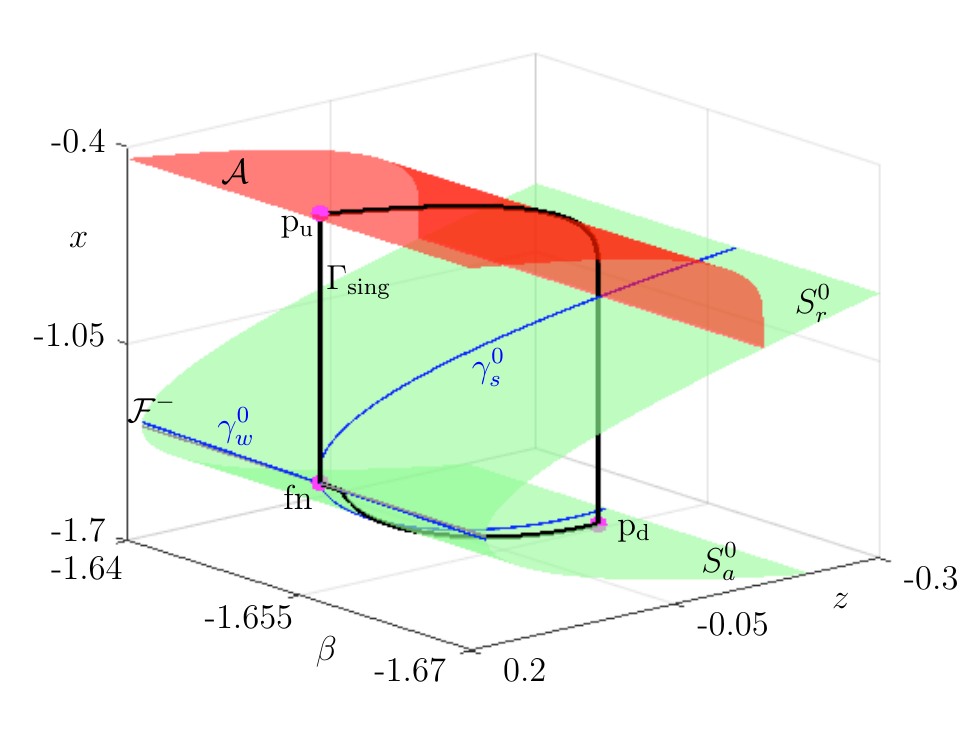}
\caption{Singular periodic orbit from which an MMBO perturbs for small enough $\eps>0$. Together with the attracting $S^0_a$ and repelling $S^0_r$ sheets of the critical manifold, the lower fold curve $\mathcal{F}^-$, and the folded node fn, also shown on this figure are the average slow nullsurface $\mathcal{A}$, the singular strong $\gamma_s^0$ and weak $\gamma_w^0$ canards, the landing-up point $\mathrm{p_u}$, the landing-down point $\mathrm{p_d}$ and the singular periodic orbit $\Gamma_{\mathrm{sing}}$.}
\label{fig:singorb}
\end{figure}

Hence the singular orbit is formed by the following segments (see Figure~\ref{fig:singorb}):
\begin{enumerate}
\item A critical fiber connecting the folded node to the landing-up point $\mathrm{p_u}$;
\item A trajectory of the averaged slow flow ending along the line of bifurcation points of the fast subsystem ending the burst;
\item A fast fiber connecting that point to landing-down point $\mathrm{p_d}$;
\item A segment of the slow flow connecting $\mathrm{p_d}$ to the folded node.
\end{enumerate}
%

\subsection{Cyclic folded-node case}
%
\begin{figure}[!t]
\centering
\includegraphics{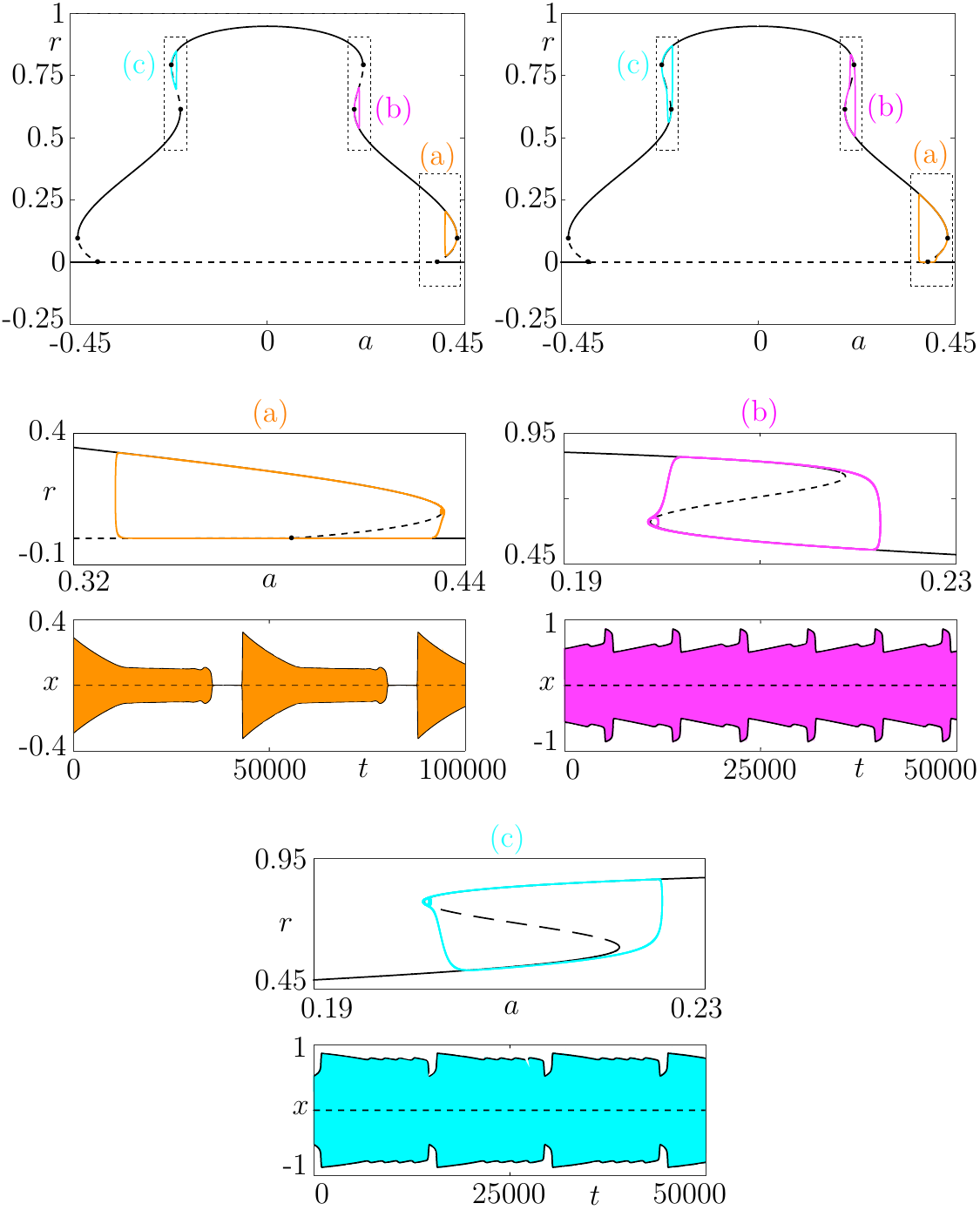}
\caption{Cyclic folded-node bursting cases. We use a polar-coordinates formulation in order to construct idealized models. The top panels show the slow-fast dissection for the amplitude $r$ of the underlying bursting model, with three different torus canard scenarios (a), (b) and (c). Adding a slow dynamics on parameter $a$ yields associated cyclic folded-node bursting scenarios for which we show both the slow-fast dissection in the $(a,r)$ plane and the $x$ time series : (a) initiated by a subcritical Hopf bifurcation; (b) terminated by a fold of cycles; (c) initiated by a fold of cycles.}
\label{fig:fncy}
\end{figure}
\begin{figure}[!t]
\centering
\includegraphics{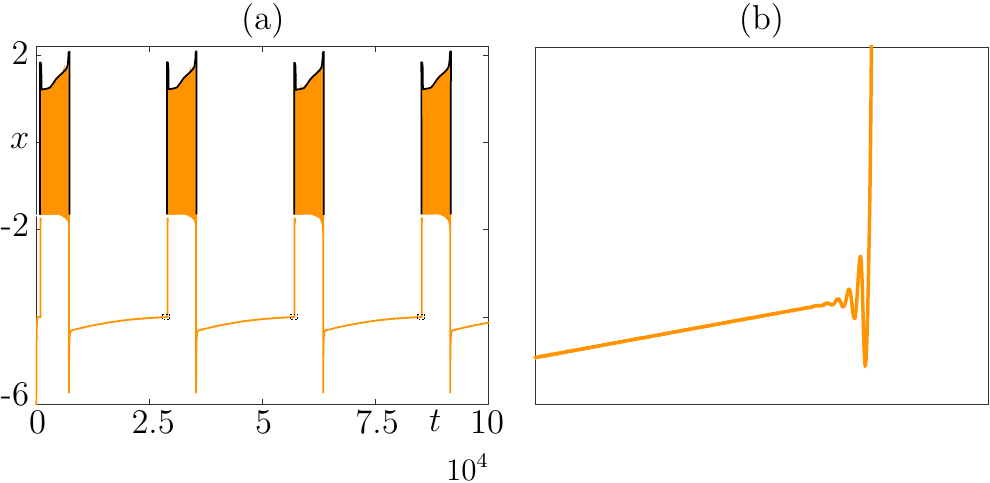}
\caption{An example of a folded-node and cyclic folded-node bursting case. The $x$-times series of the folded-node / cyclic folded-node bursting solution is shown; panel (b) is a zoom of panel (a) near the folded node.}
\label{fig:fnfncy}
\end{figure}

\mdsr{In the same spirit as in Section~\ref{sec:classicalFN}, one can construct interesting bursting rhythms where the slow oscillations occur on the envelope of the burst and this is due to what we will \red{denote} \textit{cyclic folded node}. Parallel to the construction of a folded-node burster system, one can construct a cyclic-folded-node burster system by considering a three-dimensional slow-fast system which possesses \textit{torus canard} solutions. Loosely speaking, this corresponds to a canard phenomenon with a fast rotation. Already mentioned by Izhikevich in~\cite{izhikevich01} in a canonical model, it \red{was later} found in a biophysical model of Cerebellar Purkinje cell exhibiting fold/fold cycle bursting~\cite{kramer08}, and subsequently analysed with more mathematical details in, e.g.,~\cite{benes11,burke12}. Even though to date not all elements of torus canard transitions have been \red{mathematically} unravelled, one can summarise this phenomenon by \red{emphasising that its key feature} corresponds to a canard explosion within a fast oscillatory motion. Instead of following slowly a family of equilibria past a fold bifurcation, the fast-oscillating system follows slowly a family of equilibria past a cyclic fold bifurcation. \red{Moreover, one can draw} a parallel between classical canards and torus canards in their role of transitional regime in neuronal dynamics: classical canards can explain the rapid transition from rest to the spiking regime, likewise torus canards can explain the rapid transition from the spiking to the bursting regime. \red{Furthermore}, torus canards are also not robust and only exist within exponentially thin parameter regions. \red{Thus, the} very same idea that leads from canard point to folded singularities, can lead from torus canard to cyclic folded-node canards, when adding a second slow variable. \red{In this way}, a cyclic folded-node can be robust even if the torus canards are not robust. This has been proposed very recently by Vo and collaborators~\cite{vo17,vo16} \red{via a specific example that links} the resulting dynamics to the amplitude-modulated bursting already mentioned in~\cite{izhikevich01,kramer08}; see also~\cite{han18} for other examples of amplitude-modulated bursting. \red{In summary, we herein} propose a taxonomy of cyclic folded-node bursting patterns, with several numerical examples, which completes our extension of the previous bursting classifications. \red{We complement this} with few examples of idealized models displaying cyclic-folded-node bursting. \red{This is achieved by considering systems expressed} in polar form, in which case the condition for cyclic folded node and then for cyclic folded-node bursting reduce to folded-node conditions on $r$; see Figure~\ref{fig:fncy}. In general, \red{it is possible to reduce} the system locally near the cyclic fold bifurcation of the fast subsystem \red{enabling the computation of} normal form coefficients (see~\cite{roberts15,roberts17,vo16,vo17}) \red{that effectively characterise} the cyclic folded-node. However, the bursting conditions have not been established in general. \red{Finally, for sake of completeness, we construct a limiting case of a non-trivial system that displays both classical folded-node bursting and cyclic-folded-node bursting, as depicted in} Figure~\ref{fig:fnfncy}.}

%

\section{Conclusion and perspectives}
\label{sec:conclusion}
\mdsr{The research topic of mathematical classification of bursting patterns was initiated with seminal papers published in the mid 1980s with three proposed classes of bursting oscillations~\cite{rinzel85,rinzel86,rinzel87}.} The key idea of comparing the fast subsystem's bifurcation diagram and the full systems' dynamics may seem natural with \red{hindsight, but in fact it was a genuine breakthrough, which shaped the way bursting oscillations have been modelled and dissected ever since}. The present manuscript follows \mdsr{these} footsteps, as well as those of the subsequent contributors on this topic~\cite{bertram95,golubitsky01,izhikevich00}\mdsr{, hence it was important to review in detail \red{these results since they form one pilar of Mathematical Neuroscience but also have impact in other fields. We then take a step forward by proposing an extension of the classification scheme, which allows to cover more types of burster systems, in particular those with two slow variables, namely folded-node bursters}. What we propose is} a conceptual framework which does not pretend to full mathematical rigour \red{(i.e. theorem derivations since its beyond the focus of the paper)} but rather it aims to \red{provide key insights that will enable} further theoretical and modelling work on the vast question of bursting. \mdsr{The extended bursting classification crucially focuses on the dynamics \red{during the} silent phase where the termination \red{of the trajectory profile} is not just a simple rise over the fold of the critical manifold but can involve subthreshold oscillations. We emphasize the importance of the slow flow in slow-fast systems with (at least) two slow variables, which was somehow previously overlooked in the context of bursting. In such two-slow-variable bursting systems, the silent phase termination is due to the presence of folded node. This scenario is known to give rise to canard solutions that organise, upon parameter variation but also, transiently, upon change of initial conditions, the number of subthreshold oscillations. This slow cycle-adding phenomenon is indeed entirely due to canards \red{and it} controls the profile of the underlying bursting oscillations. \red{Importantly,} it does so in a robust manner in the sense that such bursting patterns with subthreshold oscillations exist over order one ranges of parameter values. Therefore, we have added to the state-of-the-art classifications of bursting patterns the cases of classical and cycle folded-node bursters, which we can summarise in~\ref{tab:fnburst} below. We also propose a few more cases to the existing classifications (refer to Appendix~\ref{sec:newfast}), which to the best of our knowledge, have not been reported before.}
\begin{table}[!h]
\begin{center}
\begin{tabular}{|l|c|c|c|}\hline
\diaghead{\theadfont Diag ColumnmnHead II}%
  {Initiation\\of the burst}{Termination\\of the burst}&
\thead{codim. 1 bif.\\of cycles}&\thead{\blue{cyclic \textbf{fn}}}&\thead{\blue{classical \textbf{fn}}}\\ \hline
\blue{classical \textbf{fn}} & \blue{\cmark},~\cite{desroches13a} & \blue{\cmark} &  \cite{vbw}+1 slow var.\\    \hline
codim. 1 bif. equilibria &  \cite{izhikevich00,rinzel85} & \blue{\cmark},~\cite{vo16} & \cite{vbw} \\    \hline
codim. 1 bif. cycles &  \cite{izhikevich00,rinzel85} & \blue{\cmark} & \cite{vbw} \\    \hline
\blue{cyclic      \textbf{fn}} & \blue{\cmark} & \texttt{modify Fig.~\ref{fig:fncy} (b)} & \blue{\cmark} \\    \hline
\end{tabular}
\vspace*{0.2cm}
\caption{Extended classification of bursting patterns.}
\label{tab:fnburst}
\end{center}
\end{table}

\mdsr{Where do we go from here? Following this initial framework for folded-node bursting, it will be important to develop this approach in the context of biophysical excitable cell models with more than one slow processes. To this extent, a very interesting question for follow-up work is to rethink \red{about} folded-node bursting dynamics from a \red{biophysical} modelling viewpoint. In all our \mdsr{idealized} models of folded-node bursting, we have added feedback terms in the second slow differential equation with both positive and negative coefficients, which tends to indicate that both positive and negative feedback loops are useful to produce the desired output behaviour. \red{In this context, we hight two interesting aspects associated with the experimental time-series that we attempted to model with our idealized model (folded-node bursting) reproduced in} Figs.~\ref{fig:MMBO_exp} and~\ref{fig:fnhom}. First, the subthreshold oscillations \red{appear to be following} the excitability threshold, which may be harder to obtain in a three-dimensional model, even though some elliptic bursting models --e.g. FitzHugh-Rinzel, Morris-Lecar as well as some MMO models-- could potentially reproduce this aspect. However, \red{these elliptic bursting models cannot} capture the second aspect. Note that our example of folded-node bursting has 3 time scales; this was done for convenience in the construction and may not be absolutely necessary. Second, the burst phase is located on a plateau (in terms of \red{neuronal} membrane potential values) compared to the quiescent phase, which is reminiscent of a square-wave type bursting. Indeed our \mdsr{idealized} folded-node bursting model \red{reproduces quite well these data and in fact it can effectively be designated as a} folded-node homoclinic bursting model. Three-dimensional elliptic bursting models, or MMO models, would not be able to capture this aspect. \red{One interesting possibility to find biophysical models with folded-node bursting dynamics is perhaps via existing models} of thalamic bursting, or \red{alternatively to extend these models to explain the observational data published in}~\cite{roy84}. In terms of application to neural dynamics, it is legitimate to ask about neural coding and the implications of folded-node dynamics within a bursting regime. There, one would want to compare spike-adding to folded-node cycle-adding. The cycle adding can quantize the slow phase duration which might have significant effect on silent phase (and therefore on active phase) durations. On the other hand, spike-adding has less impact on macroscopic timing and less impact if a spike is added to a burst of several, say, 6 or more, spikes. A single spike added in a 2-4 spike burst might have coding contributions (synaptic transmission) but less so if there are already more than 6 spikes in a burst. These questions go beyond the scope of the present paper but are certainly of direct interest for follow-up work. \red{On the theoretical side, as aforementioned} we do not claim to have reached mathematical rigour but rather to have introduced a new framework for analysing bursting oscillations with two slow variables. There are clearly several open avenues for \red{rigorous theorem-driven directions} to be explored. \red{For instance,} we have \red{elucidated} how to construct a singular orbit that will perturb to a folded-node bursting orbit, but, proving this perturbation result is not immediate. In general, proving rigorously the existence of canards in this 4D context and how both 3D parts (subthreshold and superthreshold) combine to organise the global dynamics require \red{non-trivial} mathematical analysis.}

\mdsr{The question of noise is also a natural one to consider. If small to moderate noise is added to a folded-node bursting systems, \red{it is likely that noise will not affect significantly the burst phase. However, it is expected that the phase} of spiking oscillations during the burst will be affected, but not the qualitative dynamics. Folded-node dynamics is known to be robust to noise, its time course is parametrically robust and noise-tolerant. The canard phenomenon accounts for subtle dynamic features like cycle-adding however the subthreshold oscillations near a folded node are robust. The noise will affect these subthreshold oscillations by modifying the rotation sector in which the trajectory falls into from one passage to the the next, however the oscillations will remain. To quantify this variability of the \red{sector of a} folded-node burster with noise, one could use results by Berglund et al.~\cite{berglund12}. However, here as well the qualitative dynamics and the key role of the slow subsystem and its folded node will remain. A rigorous understanding of the impact of noise on a folded-node burster model is certainly an interesting question that goes beyond the scope of the present work.}

\mdsr{Finally, the question of bursting dynamics with at least two slow variables and more than two timescales is also of interest and related to the present work. As \red{aforementioned}, in the limit of folded-saddle-node singularities, small subthreshold oscillations will remain and increase in number and shape. In the context of slow-fast systems with two slow variables, this scenario is well-known to be akin to three-timescale dynamics~\cite{krupa10}. The associated bifurcation structure is already involved in the three-dimensional setup, with involvement of adding organizing centers such as singular Hopf bifurcation points~\cite{guckenheimer08}. \red{Thus,} it is to be expected that the folded-saddle-node bursting profiles will be more rich and \red{complex} to fully describe than the folded-node bursting cases presented \red{herein}. Yet, the underlying robust mechanism that gives a bursting pattern and requires the analysis of both slow and fast subsystem will be similar as the one proposed in the present work. A full analysis of this limiting case is a very interesting and natural question for future work. Besides, bursting systems with more than two timescales have recently gained further interest in link with canard solutions~\cite{desroches18,krupa12,letson17,nan15}, where the additional timescales bring more structure to the system and allow for further geometric singular perturbation analysis. Such approaches would certainly shed further light onto folded-node bursting dynamics as presented here and we regard it as a natural and interesting question for future work.}

\appendix
\section{Novel bursting cases in the state-of-the-art classification scheme}
\label{sec:newfast}
\begin{figure}[!t]
\centering
\includegraphics{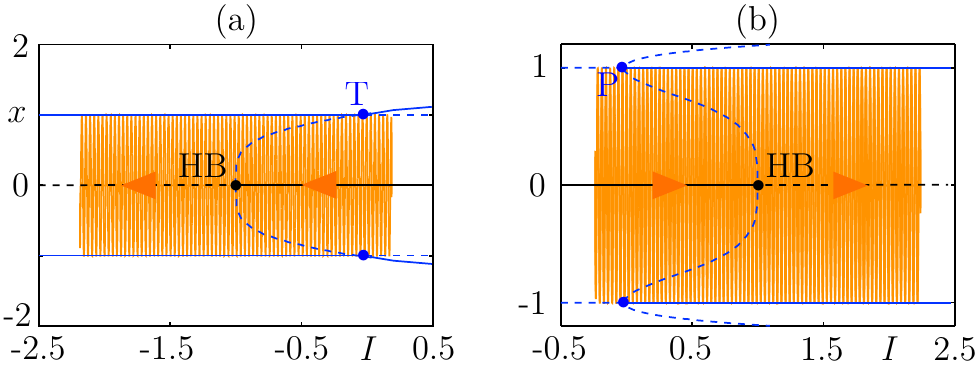}
\caption{(a1)-(b1): ``SubHopf/Transcritical of cycles'' bursting in the system: $\dot{x}=-y-x(1-r^2)(1+I-r^2)$, $\dot{y}=x-y(1-r^2)(1+I-r^2)$, $\dot{I}=\eps(r-\alpha)$, with $r=\sqrt{x^2+y^2}$ and $\alpha=0.5$. (a2)-(b2): ``SubHopf/Pitchfork of cycles'' bursting in the system: $\dot{x}=-y+x(1-r^2)(I-(r^2-1)^2)$, $\dot{y}=x+y(1-r^2)(I-(r^2-1)^2)$, $\dot{I}=\eps(\alpha-r)$, with $r=\sqrt{x^2+y^2}$ and $\alpha=0.5$.}
\label{fig:transpitch}
\end{figure}
%
Noteworthy, researchers have been challenged over the recent years to provide an exhaustive list of models (both \mdsr{idealized} and biophysical) that complete the bursting classification list predicted by the state-of-the-art classifications. This challenge is driven by the wealth of novel experimental data showing time-series predicted by the state-of-the-art classifications but for which no model has yet been implemented. Thus before moving towards  the new classification framework for bursting oscillations, we asked the question if there are bursting scenarios predicted by the Rinzel-Izhikevich and Bertram-Golubitsky classification system but which have not yet been implemented (i.e. to the best of our knowledge uncharted within the literature). Indeed, we identified a few cases in the literature and also constructed our own novel \mdsr{idealized} models but by no means this implies we have finally exhausted all cases predicted by state-of-the-art classifications. Specifically, the examples that we have managed to construct involve systems that have either one slow variable or two slow variables, which are discussed subsequently.

\subsection{Novel bursting cases with one slow variable}
%
\begin{figure}[!t]
\centering
\includegraphics{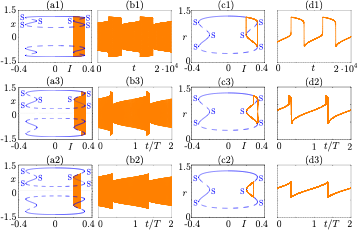}
\caption{The constructed system has the following equations: $r'=2r(0.01-(r-c)^4+\alpha^2((r-b)^2-I^2))$, $I'=\eps(r-k)$ with the following parameter values: $c=0.8$, $\alpha=0.5$, $b=0.8$, $\eps=10^{-4}$ and varying $k$ for different configurations. (a1)-(d1): ($k=0.9$) ``Fold of cycles/fold of cycles'' Isola bursting. The underlying fast subsystem bifurcation branch lies on an isola. (a2)-(d2): ($k\approx0.79996867364$) Headless torus canard on the limit cycle branch leading to isola bursting orbit. (a3)-(d3): ($k\approx0.79996867364$) Torus canard with head on the limit cycle branch leading to isola bursting orbit. Lower panels show the envelope of the bursting cycle plotted on the upper panels, with $r^2=x^2+y^2$.}
\label{fig:isola}
\end{figure}
We begin by listing all the cases we have identified in the literature as recent novel bursting models that possess one main slow variable and, hence, rely upon bistability of the fast subsystem. Recently, a biophysical $\beta$-cell model was proposed by Yildirim et al.~\cite{yildirim17}, where for the first time they were able to present the construction of a model whereby the second bifurcation (of the fast subsystem) that controls the termination of the burst is a period-doubling bifurcation. Also recently, Han et al.~\cite{han14} derived a number of new cases in which the bursting oscillations is initiated via a transcritical bifurcation of equilibria or through a pitchfork bifurcation of equilibria. We subsequently propose three new cases, which to the best of our knowledge have not been reported in the literature, neither in biophysical models nor in \mdsr{idealized} models. These cases involve, \textit{transcritical bifurcation of cycles}, \textit{pitchfork bifurcation of cycles} and \textit{isola of limit cycles} in the fast subsystem mediating bursting cycles in the full 3D system. For the first two cases, we implement \mdsr{idealized} models by first considering a normal form of transcritical and pitchfork bifurcation written in terms of the radius $r$ of the corresponding cycles of the fast subsystem. To this end, we write the polar-coordinate version of the fast subsystem and then derive from it the cartesian coordinate equations for $x$ and $y$. Finally, we append to each $(x,y)$ system a simple slow equation allowing for a periodic passage through a bistable region of the fast subsystem. These new models are depicted in Fig.~\ref{fig:transpitch}, where panels (a1)-(b1) is the transcritical case and panels (a2)-(b2) correspond to the pitchfork case. As in previous figures, the left panels show the classical slow-fast dissection and the projection of the bursting cycle of the full system onto the bifurcation diagram of the fast subsystem; the right panels show the corresponding time series of this bursting cycle for the fast variable $x$.
The case for \textit{isola} bursting is an intriguing case. First recall that an \textit{isola} refers to a bifurcation branch of solutions that is closed and isolated in parameter space. These can correspond to families of equilibria or limit cycles. Isolas are important objects well described by singularity theory~\cite{golubitsky85} and which also naturally appear, in the limit-cycle case, in multiple-timescale systems with at least three variables; see for instance~\cite{desroches12}. An isola (of equilibria or of limit cycles) has at least two fold points, which correspond to saddle-node bifurcation points. However, one can imagine a more complicated scenario where the isola has additional fold points (necessarily a multiple of 2). The unusual aspect of a bursting oscillations via an isola is that the Rinzel-Izhikevich classification system classify a bursting oscillation in the same way if it is either mediated via an isola or a regular bifurcation solution branch. That is, the bifurcation points of the fast subsystem that encode the initiation and termination of the bursts would be the same, but the underlying type of solution branches would be invisible to the classification. The is a contentious point (or perhaps not so critical depending on the point of view), but there could be scenarios where from a biophysical point of view, it is crucial to consider isola as also part of the explanatory arguments to explain an underlying experimental observation.
Here we construct a system possessing an isola of limit cycles with a total of six fold points around it (see Figures~\ref{fig:isola}). The construction of the model follows the same strategy as before, whereby we start by first writing the system in terms of polar-coordinate system and only the equation for the radius $r$ is non trivial. Subsequently, we add a simple slow equation with feedback onto one of the two fast variables, which then induces bursting cycles that move, in the slow-fast dissection sense, along the isola of cycles previously described. Note, the stability along this isola changes through these fold points (see the left panels of Figures~\ref{fig:isola}). Moreover, the isola is the critical manifold of this bursting system and we choose to denote the resulting bursting pattern as \textit{isola bursting}. We re-emphasise, this bursting pattern is not entirely new \emph{per se} since it falls into the ``fold cycle / fold cycle'' case of Izhikevich's classification~\cite{izhikevich00}. However, the two fold bifurcations of cycles of the fast subsystem happen to lie on an isola and indeed this is the novelty of this scenario; see Fig.~\ref{fig:isola}. This bursting pattern can be seen as an amplitude-modulated bursting~\cite{kramer08,vo16} and it appears through a \textit{torus canard}~\cite{burke12} phenomenon when varying a parameter that controls the slow nullcline. In particular, two (non-robust) isola bursting canard scenarios occur on the way to the robust isola bursting scenario shown in Fig.~\ref{fig:isola} panels (a1)-(d1): namely, (torus) canard isola bursting without head, in Fig.~\ref{fig:isola} (a2)-(d2), and (torus) canard isola bursting with head, in Fig.~\ref{fig:isola} (a3)-(d3).

\subsection{Novel bursting cases with two slow variables}
%
\begin{figure}[!t]
\centering
\includegraphics{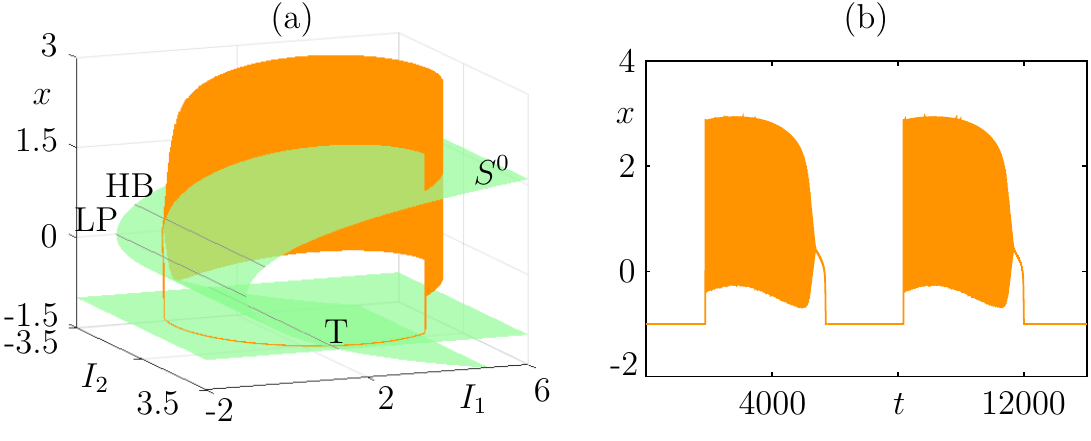}
\caption{``Transcritical of equilibria/Hopf'' bursting with two slow variables: same signature but amplitude modulation during the burst. The model we used for this simulation is: $\dot{x}=(x-k)(y-\beta x^4+3x^2+I_2)/c$, $\dot{y}=1-5x^2-y$, with: $c=4$, $\beta=0.3$, $k=-1$, $\eps=0.001$ and $(I_1,I_2)$ being a slow harmonic oscillator of speed $\eps$ centered at $(0,2)$ and of amplitude $\mu=\sqrt{10}$.}
\label{fig:trans2sl}
\end{figure}
\begin{table}[!t]
\begin{center}
\begin{tabular}{|l|c|c|c|}\hline
\diaghead{\theadfont Diag ColumnmnHead II}%
  {Initiation\\of the burst}{Termination\\of the burst}&
\thead{\blue{pitchfork}\\\blue{cycles}}&\thead{\blue{transcritical}\\\blue{cycles}}&\thead{\blue{fold of cycles}\\\blue{on isola}}\\ \hline
subcritical Hopf & \blue{\cmark} & \blue{\cmark} &  \\    \hline
\blue{fold of cycles on isola} & & & \blue{\cmark} \\    \hline
\end{tabular}
\caption{New bursting cases within the state-of-the-art classification.}
\label{tab:newinold}
\end{center}
\end{table}
%
We now turn attention to bursting models with two slow variables, but that nevertheless fall under the state-of-the art classification. To the best of our knowledge we find in the literature two principle mechanisms. The first case is the so-called parabolic bursting, which was in fact discovered by Rinzel~\cite{rinzel86,rinzel87}. This is a special case whereby that both the initiation and termination of the bursting oscillation is mediated by SNIC bifurcations of the fast subsystem. The fast subsystem does not possess bistability and the underlying bursting oscillations are not driven by hysteresis but rather the two slow variables produce a slow oscillations between themselves, which then drives the full system from one side of the curve of SNIC bifurcation points to the other. The second scenario considered in the literature is the more intuitive case, whereby the fast subsystem is bistable~\cite{smolen93}. This is in fact similar to the case of systems with only a single slow variable, however, in this model the second variable does not introduce further complex dynamics. Rather, the bifurcations of the fast subsystem are extended (via the second slow variable) to a line of bifurcations. Specifically, in the model developed in~\cite{smolen93} two cases are presented. One case illustrates how a bursting cycle is initiated via a line of fold bifurcation points and terminates via a line of homoclinic bifurcation points. In the alternative example, the bursting cycle also initiates via a line of fold bifurcation points but terminates through a family of Hopf bifurcations. We further propose a novel bursting scenario with two slow variables that also falls under the state-of-the art classification. The first point to notice is that could potentially consider a second slow variable that is capable of slow modulating the fast dynamics during the burst. This strategy can potentially prove useful in some modelling contexts where a hidden slow modes may act upon the fast dynamics by modulating its amplitude. To show this concept we construct an example whereby the the fast subsystem is bistable and the bifurcation of the fast subsystem that initiates the burst is a line of transcritical bifurcation of equilbria. Thus a slow passage through it, will cause a delayed exchange of stability hence containing a canard segment. Moreover, a line of Hopf bifurcation points terminates the burst. This models construction follows very closely the concepts developed in~\cite{han14}, however therein they only considered one slow variable. In contrast we consider a second slow variable and prescribe a dynamics (in our case harmonic oscillations), which then induces a slow modulation of the amplitude of the bursting oscillations (see  Fig.~\ref{fig:trans2sl}). Similar behaviour can be found in three-dimensional bursters (e.g. square-wave or elliptic) where one considers an additional slow process so that both slow variables oscillate together.
We close this section by summarising in Table~\ref{tab:newinold} the new cases that we have constructed. In blue we indicate the bifurcation of the bursting pattern that is novel. Note that for the second bifurcation of the fast subsystem we chose to implement a subcritical Hopf bifurcation, but it is to be expected that similar bursting patterns can be found and/or constructed with any other codimension-one bifurcation instead. Finally, we do not include within this table our example of two-slow-variable burster (or variants of it) since the second slow variable does not change the signature of the bursting but only allows for a slow amplitude-modulation during the burst.

\section*{Acknowledgements}
SR would like to acknowledge Ikerbasque (The Basque Foundation for Science) and moreover, this
research is supported by the Basque Government through the BERC 2018-2021 program and by the
Spanish State Research Agency through BCAM Severo Ochoa excellence accreditation SEV-2017-0718
and through project RTI2018-093860-B-C21 funded by (AEI/FEDER, UE) and acronym ``MathNEURO''.

\bibliography{references}
\bibliographystyle{plain}

\end{document}